\RequirePackage{etex}
\documentclass{article}
\usepackage[utf8]{inputenc}
\usepackage{inputenc}
\usepackage{amsmath}
\usepackage{amssymb}
\usepackage{ upgreek }
\usepackage{graphicx}
\usepackage{layout}
\usepackage{dsfont}
\usepackage{bbold}
\usepackage{xspace}
\usepackage[english]{babel}
\usepackage[top=2cm, bottom=3cm, left=2cm, right=2cm]{geometry}
\usepackage{multicol}
\usepackage{amsthm}
\usepackage{mathrsfs}
\usepackage{tikz}
\usepackage{tikz-network}
\usepackage[toc,page]{appendix}
\usepackage{todonotes}
\usepackage{complexity}
\usepackage{comment}
\usepackage{fullpage}
\usepackage{subcaption}
\usepackage{amsmath,amsfonts, amssymb, bm}
\usepackage{enumitem}
\usepackage[hidelinks]{hyperref}
\usepackage{ulem}
\usepackage{authblk}

\def\rotateclockwise#1{
  \newdimen\xrw
  \pgfextractx{\xrw}{#1}
  \newdimen\yrw
  \pgfextracty{\yrw}{#1}
  \pgfpoint{\yrw}{-\xrw}
}

\def\rotatecounterclockwise#1{
  \newdimen\xrcw
  \pgfextractx{\xrcw}{#1}
  \newdimen\yrcw
  \pgfextracty{\yrcw}{#1}
  \pgfpoint{-\yrcw}{\xrcw}
}

\def\outsidespacerpgfclockwise#1#2#3{
  \pgfpointscale{#3}{
    \rotateclockwise{
      \pgfpointnormalised{
        \pgfpointdiff{#1}{#2}}}}
}

\def\outsidespacerpgfcounterclockwise#1#2#3{
  \pgfpointscale{#3}{
    \rotatecounterclockwise{
      \pgfpointnormalised{
        \pgfpointdiff{#1}{#2}}}}
}

\def\outsidepgfclockwise#1#2#3{
  \pgfpointadd{#2}{\outsidespacerpgfclockwise{#1}{#2}{#3}}
}

\def\outsidepgfcounterclockwise#1#2#3{
  \pgfpointadd{#2}{\outsidespacerpgfcounterclockwise{#1}{#2}{#3}}
}

\def\outside#1#2#3{
  ($ (#2) ! #3 ! -90 : (#1) $)
}

\def\cornerpgf#1#2#3#4{
  \pgfextra{
    \pgfmathanglebetweenpoints{#2}{\outsidepgfcounterclockwise{#1}{#2}{#4}}
    \let\anglea\pgfmathresult
    \let\startangle\pgfmathresult

    \pgfmathanglebetweenpoints{#2}{\outsidepgfclockwise{#3}{#2}{#4}}
    \pgfmathparse{\pgfmathresult - \anglea}
    \pgfmathroundto{\pgfmathresult}
    \let\arcangle\pgfmathresult
    \ifthenelse{180=\arcangle \or 180<\arcangle}{
      \pgfmathparse{-360 + \arcangle}}{
      \pgfmathparse{\arcangle}}
    \let\deltaangle\pgfmathresult

    \newdimen\x
    \pgfextractx{\x}{\outsidepgfcounterclockwise{#1}{#2}{#4}}
    \newdimen\y
    \pgfextracty{\y}{\outsidepgfcounterclockwise{#1}{#2}{#4}}
  }
  -- (\x,\y) arc [start angle=\startangle, delta angle=\deltaangle, radius=#4]
}

\def\corner#1#2#3#4{
  \cornerpgf{\pgfpointanchor{#1}{center}}{\pgfpointanchor{#2}{center}}{\pgfpointanchor{#3}{center}}{#4}
}

\def\hedgeiii#1#2#3#4{
  \outside{#1}{#2}{#4} \corner{#1}{#2}{#3}{#4} \corner{#2}{#3}{#1}{#4} \corner{#3}{#1}{#2}{#4} -- cycle
}

\def\hedgeiiii#1#2#3#4#5{
  \outside{#1}{#2}{#5} \corner{#1}{#2}{#3}{#5} \corner{#2}{#3}{#4}{#5} \corner{#3}{#4}{#1}{#5} \corner{#4}{#1}{#2}{#5} -- cycle
}

\def\hedgem#1#2#3#4{
  
  \outside{#1}{#2}{#4}
  \pgfextra{
    \def\hgnodea{#1}
    \def\hgnodeb{#2}
  }
  foreach \c in {#3} {
    \corner{\hgnodea}{\hgnodeb}{\c}{#4}
    \pgfextra{
      \global\let\hgnodea\hgnodeb
      \global\let\hgnodeb\c
    }
  }
  \corner{\hgnodea}{\hgnodeb}{#1}{#4}
  \corner{\hgnodeb}{#1}{#2}{#4}
  -- cycle
}

\def\hedgeii#1#2#3{
  \hedgem{#1}{#2}{}{#3}
}

\newtheorem{theorem}{Theorem}[section]
\newtheorem{lemma}[theorem]{Lemma}
\newtheorem{claim}[theorem]{Claim}
\newtheorem{proposition}[theorem]{Proposition}
\newtheorem{observation}[theorem]{Observation}
\newtheorem{corollary}[theorem]{Corollary}

\newtheorem{problem}[theorem]{Problem}

\newcommand{\etal}{{\em et al.}\xspace}

\newcommand{\Input}{\textbf{Input}}
\newcommand{\Parameter}{\textbf{Parameter}}
\newcommand{\Question}{\textbf{Question}}

\newcommand{\qedclaim}{\hfill $\diamond$ \medskip}

\newenvironment{proofclaimitem}{\noindent{\em Proof of the claim.}}{}

\newcommand{\MBDG}{{\sc MB Dom Game}\xspace}

\newcommand{\oD}{\mathcal{D}}
\newcommand{\oN}{\mathcal{N}}
\newcommand{\oS}{\mathcal{S}}
\newcommand{\MBW}{\forall^{\neq}\text{-}\Sigma_2}
\newcommand{\MBWU}{\forall^{\neq}\text{-}\Sigma_{2,u}}
\newcommand{\MBWun}{\forall^{\neq}\text{-}\Sigma_{2,1}}
\newcommand{\shortMBM}{{\sc Short$_{\text{Maker}}$ MB Pos Game}\xspace}
\newcommand{\shortMBB}{{\sc Short$_{\text{Breaker}}$ MB Pos Game}\xspace}
\newcommand{\shortDD}{{\sc Short$_{\text{Dominator}}$ MB Dom Game}\xspace}
\newcommand{\shortDS}{{\sc Short$_{\text{Staller}}$ MB Dom Game}\xspace}
\newcommand{\kdomset}{{\sc $k$-Dominating Set}\xspace}
\newcommand{\hxf}{\hyp = (\som, \WS)}
\newcommand{\hyp}{\mathcal{H}}
\newcommand{\hypdom}{\mathcal{H}_{\text{dom}}}
\newcommand{\hypneighbor}{\mathcal{H}_{N}}
\newcommand{\WS}{\mathcal{F}}
\newcommand{\som}{\mathcal{X}}
\newcommand{\alphabet}{\mathcal{A}}
\newcommand{\join}{\bowtie}
\newcommand{\oDom}{\mathcal{D}}
\newcommand{\strat}{\mathcal{S}}

\tikzstyle{v}=[circle,inner sep=0, minimum size =6 pt, line width = 1pt, draw=black, fill=black, text= white]
\tikzstyle{inv2}=[circle,inner sep=0, minimum size = 18 pt, line width = 0pt, draw=white, fill=white, text= black]
\tikzstyle{R}=[circle,inner sep=0, minimum size =7 pt, line width = 1 pt, draw=red, fill = red]
\tikzstyle{inv}=[circle,inner sep=0, minimum size = 0 pt, line width = 0pt, draw=white, fill=white, text= black]

\title{On the parameterized complexity of the Maker-Breaker domination game
\footnote{This research was supported by the ANR project P-GASE (ANR-21-CE48-0001-01) and Kempe Foundation Grant No. JCSMK24-515 (Sweden).}
}
\author[1]{Guillaume Bagan}
\author[2]{Mathieu Hilaire}
\author[3]{Nacim Oijid}
\author[1]{Aline Parreau}

\affil[1]{Univ Lyon, CNRS, UCBL, INSA Lyon, LIRIS, UMR5205, F-69622 Villeurbanne, France}
\affil[2]{Univ. Bordeaux, Bordeaux INP, LaBRI UMR CNRS 5800, F-33400, Talence, France.}
\affil[3]{Umeå University, Sweden}

\date{}

\begin{document}

\maketitle

\begin{abstract}

Since its introduction as a Maker-Breaker positional game by Duchêne \etal in 2020, the Maker-Breaker domination game has become one of the most studied positional games on vertices. In this game, two players, Dominator and Staller, alternately claim an unclaimed vertex of a given graph $G$. If at some point the set of vertices claimed by Dominator is a dominating set, she wins; otherwise, i.e. if Staller manages to isolate a vertex by claiming all its closed neighborhood, Staller wins.  

Given a graph $G$ and a first player, Dominator or Staller must have a winning strategy. We are interested in the computational complexity of determining which player has such a strategy.
This problem is known to be \PSPACE-complete on bipartite graphs of bounded degree and split graphs; polynomial on cographs, outerplanar graphs, and block graphs; and in \NP\ for interval graphs.

In this paper, we consider the parameterized complexity of this game. We start by considering as a parameter the number of moves of both players. 
We prove that for the general framework of Maker-Breaker positional games in hypergraphs, determining whether Breaker can claim a transversal of the hypergraph in $k$ moves is \W[2]-complete, in contrast to the problem of determining whether Maker can claim all the vertices of a hyperedge in $k$ moves, which is known to be \W[1]-complete since 2017. 
These two hardness results are then applied to the Maker-Breaker domination game, proving that it is \W[2]-complete to decide if Dominator can dominate the graph in $k$ moves and \W[1]-complete to decide if Staller can isolate a vertex in $k$ moves. 
Next, 
we provide \FPT\ algorithms for the Maker-Breaker domination game parameterized by the neighborhood diversity, the modular width, the $P_4$-fewness, the distance to cluster, and the feedback edge number.
\end{abstract}

\section{Introduction}

The Maker-Breaker domination game \cite{makerbreaker} is played on a graph $G$. Two players, Dominator and Staller, are alternately claiming vertices of the graph. Dominator wins if at some point her vertices form a dominating set, otherwise, i.e. if Staller manages to have a vertex and all its neighbors, Staller wins. This game belongs to the family of {\em Maker-Breaker positional games}, introduced in 1963 by Hales and Jewett \cite{Hales1963} and independently in 1973 by Erd\"os and Selfridge \cite{erdos}. A Maker-Breaker positional game is played on a hypergraph $\mathcal H$. The two players, Maker and Breaker, alternatively claimed vertices of $\mathcal H$. Maker wins if at some point she has claimed all the vertices of a hyperedge of $\mathcal H$. Otherwise, Breaker wins. To model the Maker-Breaker domination game played on a graph $G$ as a positional game on $\mathcal H$, the vertices of $\mathcal H$ are the vertices of $G$ and the hyperedges can be either the dominating sets of $G$ (then Dominator plays the role of Maker) or can be the closed neighborhoods of $G$ (and then Dominator plays the role of Breaker).
 
 As two-player finite games with perfect information, no chance and no draw, in a given Maker-Breaker positional game, one of the players has a winning strategy.
Thus, one can consider the following decision problem \MBDG: given a graph $G$ and the first player (Dominator or Staller) as input, who has a winning strategy in the Maker-Breaker domination game played on $G$?
This problem has been proved to be \PSPACE-complete even when restricted to bipartite or chordal graphs in the paper introducing this game~\cite{makerbreaker}, and later restricted to bounded degree graphs~\cite{Oijid2025}. On the positive side, it has been proved to be polynomial for forests and cographs~\cite{makerbreaker} and later for block graphs, outerplanar graphs, regular graphs and unit interval graphs~\cite{Intervalle}. For interval graphs, \MBDG has been shown to be in \NP\ but it is an open question to know if \MBDG is actually polynomial or not~\cite{Intervalle}.

In this paper, we consider the parameterized complexity of \MBDG: which parameters make the problem hard or easy? The most natural parameter to consider in the context of positional games is the number of moves. The decision problem becomes the following: is it possible to win for a given player in at most $k$ moves? For Maker-Breaker positional games,  Bonnet et al. \cite{Bonnet2017} have proved using a new fragment of first-order logic that deciding if Maker can win in at most $k$ moves is \W[1]-complete. Some results for particular games are also known: {\sc Generalized Hex} is also \W[1]-complete~\cite{Bonnet2017} but {\sc Hex}, {\sc short $k$-connect} and the star-game are \FPT\ \cite{Bonnet2017, Bonnet2016, edgegame} when parameterized by the number of moves of Maker. 
It seems that nothing is known when the parameter is the number of moves of Breaker, i.e if Breaker has to claim a vertex in each hyperedge in a limited number of moves.
Considering only structural parameters, the size of the hyperedges has been studied. Deciding the winner of a Maker-Breaker positional game with hyperedges of size 3 is polynomial~\cite{galliot} but \PSPACE-complete already for hyperedges of size 4~\cite{galliot2025} and thus one cannot hope for \FPT\ results for this parameter. 
For the specific domination game, an XP algorithm is known for $k$-nested interval graphs~\cite{Intervalle}. 
 In this paper, we consider classical graph parameters for \FPT\ algorithms in algorithmic graph theory, namely modular-width~\cite{Lampis2012, HATANAKA2018}, $P_4$-fewness~\cite{babel2001efficient, CAMPOS2014}, distance to cluster~\cite{Gomes2021}, and the feedback edge number~\cite{kellerhals2020}

\subsection{Results.}
In this paper, considering Staller as a Maker player of a positional game, we prove that \MBDG parameterized by the number of moves of Staller is \W[1]-complete. But it turns out to be \W[2]-complete when parameterized by the number of moves of Dominator to win. To prove this result, we extend the fragment of logic given by Bonnet et al. \cite{Bonnet2017}. We actually prove a stronger result: deciding if Breaker can win in $k$ moves in a Maker-Breaker positional game is \W[2]-complete. This result is surprising in contrast to its counterpart parameterized by the number of moves of Maker, which is \W[1]-complete. This can be explained by the fact that checking whether Maker has won is a local property (a hyperedge is claimed by Maker), whereas checking whether Breaker has won is a global property (each hyperedge contains a vertex of Breaker).

Next, we consider structural graph parameters. One could naturally consider the size of a dominating set, but this parameter is not relevant since even for graphs with dominating sets of size 1 it is already \PSPACE-complete to compute the outcome. Since \MBDG is polynomial in the families of forests and cographs, we consider parameters that generalize these families. In particular, cographs can be easily solved because this game behaves well with respect to union and join operations. We emphasize this fact by proving that a module can be replaced by a graph with at most three vertices. As a corollary, we obtain a linear kernel when the problem is parameterized by the neighborhood diversity of the graph and an \FPT\ algorithm considering the modular-width as parameter. Another property of cographs is that they do not contain $P_4$ as an induced subgraph. One can generalize this property by considering the {\em $P_4$-fewness} of a graph~\cite{P4fewness}. We obtain an \FPT\ algorithm using this parameter.
We then consider the {\em distance to cluster}, i.e. the number of vertices to delete to obtain a union of cliques. If this parameter is small compared to the size of the graph, many twin vertices appear, allowing us to reduce the graph to a kernel which implies a \FPT\ algorithm for this parameter. Finally, since forests can be solved in linear time, we consider the minimum size of a feedback edge set. When this parameter is small, long induced paths appear in the graph. We prove that one can reduce induced paths to paths of length at most 7 without changing the outcome, which again lead to a kernel and an \FPT\ algorithm for the feedback edge number.

\subsection{Outline of the paper.} 
\begin{figure}[ht]
\centering
\scalebox{.8}{
\begin{tikzpicture}

\draw (0,0) node[text width=5em, text centered] (vc) {vertex cover};

\draw (0,-1.5) node[text width=5em, text centered] (fes) {feedback edge number};

\draw (0,-3) node[text width=5em, text centered] (fvs) {feedback vertex set};

\draw (0,-4.5) node[text width=5em, text centered] (tw) {treewidth};

\draw (0,-6) node[text width=5em, text centered] (cw) {twin-width};

\draw (-4,-1.5) node[text width=5em, text centered] (nd) {neighborhood diversity};

\draw (-8,-1.5) node[text width=5em, text centered] (pf) {$P_4$-fewness};

\draw (-5.5,-3) node[text width=5em, text centered] (mw) {modular width};
\draw (-2.5,-3) node[text width=5em, text centered] (dtc) {distance to cluster};

\draw (4,-1.5) node[text width=5em, text centered] (dtuos) {distance to union of stars};

\draw[style = dashed] (-9,-4) -- (-1.3,-4);
\draw[style = dashed] (-1.3,-2.2) -- (-1.3,-4);
\draw[style = dashed] (-1.3,-2.2) -- (2.3,-2.2);
\draw[style = dashed] (2.3,-.5) -- (2.3,-2.2);
\draw[style = dashed] (2.3,-.5) -- (5,-.5);

\draw (-8.5,0) ellipse (1.5cm and .7cm);
\draw  (-8.5,0) node[inv] {\Large \FPT}; 

\draw (4,-4.5) ellipse (1.5cm and .7cm);
\draw  (4,-4.5) node[inv] {\Large Open}; 

\draw[->] (vc) -- (nd);
\draw[->] (vc) -- (pf);
\draw[->] (pf) -- (-8, -6) --  (cw);
\draw[->] (vc) -- (fes);
\draw[->] (vc) -- (dtuos);
\draw[->] (nd) -- (mw);
\draw[->] (nd) -- (dtc);
\draw[->] (mw) -- (cw);
\draw[->] (dtc) -- (cw);
\draw[->] (fes) -- (fvs);
\draw[->] (fvs) -- (tw);
\draw[->] (tw) -- (cw);
\draw[->] (dtuos) -- (fvs);

\end{tikzpicture}}

\caption{Summary of our results on structural parameters. If a parameter is bounded for a class of graphs, it is bounded for its descendants in the diagram.}
\label{fig:hasse diagram}

\end{figure}
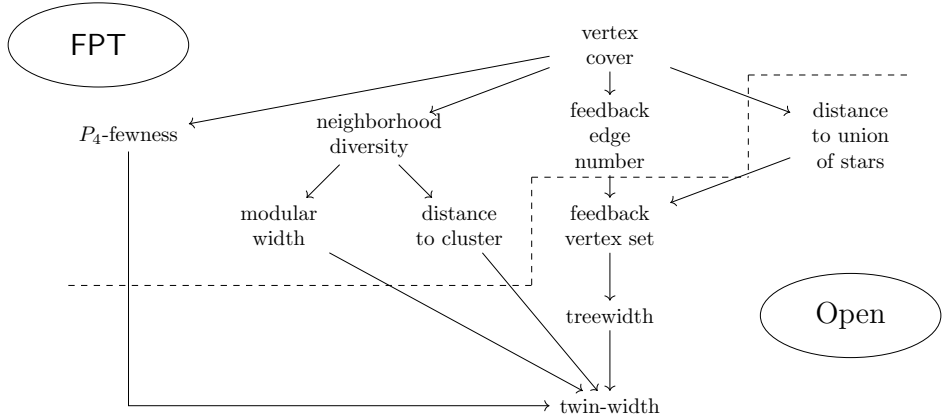

The paper is organized as follows. We start by providing formal definitions and preliminaries in Section \ref{sec:prelim}. Section~\ref{sec:moves} focuses on the number of moves. In Section \ref{sec:modules}, we consider modules. After giving a general lemma to replace modules in graphs without changing the outcome, we give FPT algorithms for the neighborhood diversity, the modular-width and the $P_4$-fewness. Section \ref{sec:cluster} and Section \ref{sec:fes} are dedicated to the parameters "Distance to Cluster" and "Feedback edge number". The diagram in Figure~\ref{fig:hasse diagram} summarizes our results on structural parameters.

\section{Preliminaries}\label{sec:prelim}

In this section, we first give formal definitions on graphs and games that we are using in the paper. We then recall some useful results on the Maker-Breaker domination game as well as technical tools we will need in some proofs.

\subsection{Graphs}

 Let $G = (V,E)$ be a simple and undirected graph with vertex set $V=V(G)$ and edge set $E=E(G)$. If $x \in V$ is a vertex, we denote by $N(x)$ the open neighborhood of $x$, and by $N[x]$ its closed neighborhood. Formally, we have $N(x) = \{y \in V | xy \in E\}$ and $N[x] = \{x\} \cup N(x)$. If $S$ is a set of vertices, we denote by $N(S)$ (resp. $N[S]$), the union of the neighborhoods of the vertices in $S$, i.e. we have $N(S) = \underset{v \in S}{\cup} N(v)$ (resp. $N[S] = \underset{v \in S}{\cup} N[v]$). A {\em leaf} is a vertex of degree 1. A {\em support vertex} is a vertex adjacent to a leaf.
 Let $S$ be a subset of vertices of $G$. The subgraph of $G$ {\em induced} by $S$, denoted by $G[S]$, is the graph with vertex set $S$ and edge set all the edges of $G$ between two vertices of $S$, i.e., the set $E\cap (S\times S)$.
 A subset $S$ of vertices is a {\em dominating set} if every vertex of $G$ is either in $S$ or adjacent to some vertex of $S$. In other words, if we have $N[S] = V$. The size of a minimum dominating set of $G$ is denoted by $\gamma(G)$. 

  A {\em module} of $G$ is a subset of vertices $M$ such that, for any $x,y \in M$, we have $N[x] \setminus M = N[y] \setminus M$. Since all the vertices of a module have the same neighborhood outside $M$, we say that a vertex $x \in V\setminus M$ is {\em adjacent to $M$}, if it is adjacent to all the vertices of $M$.
 If $u,v \in V$ are two vertices such that $\{u,v\}$ is a module, we say that $u$ and $v$ {\em have the same type} or are {\em twins}.
If $M$ is a module of $G$ and $H$ a graph, we say that $G'$ is obtained by {\em replacing} $M$ by $H$ if $G'$ is obtained from $G$ by removing the vertices of $M$, adding a copy of $H$, and adding all the edges between vertices of $H$ and vertices of $G \setminus M$ adjacent to $M$. 

Let $G$ and $H$ be two graphs. The {\em union} of $G$ and $H$, denoted $G \cup H$, is the graph obtained by taking disjoint copies of $G$ and $H$. The {\em join} of $G$ and $H$, denoted $G \join H$, is the graph obtained by taking one copy of $G$ and one copy of $H$ and adding all the edges between them. A graph is said to be a {\em cograph} if it can be constructed from isolated vertices using only union and join operations.

We denote by $K_{n}$ the complete graph on $n$ vertices and by $P_n$ the path on $n$ vertices.

\subsection{Games}

Let $\hyp=(\som,\WS)$ be a hypergraph with vertex set $\som=\som(\hyp)$ and hyperedge set $\WS=\WS(\hyp)$. A {\em transversal} of $\hyp$ is a subset of vertices that is intersecting each hyperedge of $\hyp$. The {\em Maker-Breaker positional game} on $\hyp$ is played by two players, Maker and Breaker. Maker and Breaker take turns claiming each unclaimed vertex of $\som$. Maker wins if she claims all the vertices of a hyperedge $f \in \WS$. Otherwise, i.e if Breaker manages to claim a {\em transversal} of $\hyp$, Breaker wins. 

Given a graph $G = (V, E)$, the {\em Maker-Breaker domination game} on $G$ is played by two players, Dominator and Staller. Dominator and Staller take turns claiming each unclaimed vertex of $V$. Dominator wins if she claims all the vertices of a dominating set of $G$. Otherwise, Staller wins. 
 The Maker-Breaker domination game played on a graph $G$ can be seen as a Maker-Breaker positional game using two different hypergraphs. One can consider the hypergraph $\hypdom(G)$ with vertex set $V$ and hyperedges all the minimum dominating sets of $G$: $\WS(\hypdom(G)) = \{S \subset V | S \text{ is a minimum dominating set of }G\}$. Then Dominator plays the role of Maker and Staller plays the role of Breaker. The other possibility
is to consider that Staller plays the role of Maker. Indeed, Staller wins if and only if he managed to claim a vertex and all its neighbors. The corresponding hypergraph $\hypneighbor(G)$ has again vertex set $V$ but its hyperedges are all the closed neighborhoods: $\WS(\hypneighbor(G)) = \{N[x] | x\in V\}$. 
Note that a graph might have an exponential number of minimum dominating sets compared to its number of vertices, thus $\hypdom(G)$ might be of exponential size whereas $\hypneighbor(G)$ always has at most $|V|$ hyperedges.

 Let $D$ and $S$ be disjoint subsets of  vertices. The {\em position} $(G, D, S)$ is the state of the Maker-Breaker domination game played in $G$ where Dominator has claimed the vertices of $D$ and Staller the vertices of~$S$.

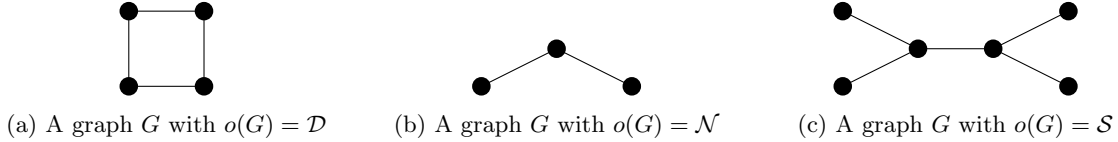
\begin{figure}[t!]
    \centering
    \begin{subfigure}[t]{0.3\textwidth}
        \centering

\begin{tikzpicture}

\coordinate (V1) at (0,0);
\coordinate (V2) at (1,0);
\coordinate (V3) at (0,1);
\coordinate (V4) at (1,1);

\draw (V1) node[v] {};
\draw (V2) node[v] {};
\draw (V3) node[v] {};
\draw (V4) node[v] {};

\draw (V1) -- (V2);
\draw (V2) -- (V4);
\draw (V3) -- (V4);
\draw (V1) -- (V3);

\end{tikzpicture}
        
        \caption{A graph $G$ with $o(G) = \oDom$}
    \end{subfigure}%
    ~ 
    \begin{subfigure}[t]{0.3\textwidth}
        \centering

\begin{tikzpicture}

\coordinate (V1) at (0,0);
\coordinate (V2) at (-1,-.5);
\coordinate (V3) at (1,-.5);

\draw (V1) node[v] {};
\draw (V2) node[v] {};
\draw (V3) node[v] {};

\draw (V1) -- (V2);
\draw (V1) -- (V3);

\end{tikzpicture}
        \caption{A graph $G$ with $o(G) = \oN$}
    \end{subfigure}
    ~
\begin{subfigure}[t]{0.3\textwidth}
        \centering

\begin{tikzpicture}

\coordinate (V1) at (0,0);
\coordinate (V2) at (-1,.5);
\coordinate (V3) at (-1,-.5);
\coordinate (V4) at (1,0);
\coordinate (V5) at (2,-.5);
\coordinate (V6) at (2,.5);

\draw (V1) node[v] {};
\draw (V2) node[v] {};
\draw (V3) node[v] {};
\draw (V4) node[v] {};
\draw (V5) node[v] {};
\draw (V6) node[v] {};

\draw (V1) -- (V2);
\draw (V1) -- (V3);
\draw (V1) -- (V4);
\draw (V4) -- (V5);
\draw (V4) -- (V6);    
\end{tikzpicture}
        
        \caption{A graph $G$ with $o(G) = \oS$}
    \end{subfigure}%
    \caption{The different outcomes in the \MBDG}
    \label{fig: outcomes}
\end{figure}

As a $2$-player perfect and finite information game, in a given position of the Maker-Breaker positional game, a player has a winning strategy: he can win no matter what moves his opponent makes. Moreover, in a Maker-Breaker positional game, a winning strategy for some player going second implies a strategy for the same player going first (see \cite{Hefetz2014} for a proof). Therefore, for any position $(G, D, S)$ of a Maker-Breaker domination game, there are three possible {\em outcomes}: either Dominator can ensure a win playing first or second (outcome $\oDom$), or Staller can ensure a win playing first or second (outcome $\oS$) or the next player to move can ensure a win (outcome $\oN$). The outcome of a position will be denoted by $o(G, D, S)$. If the position is the starting position, i.e. if $D = S = \emptyset$, we simply note the outcome by $o(G)$. See Figure~\ref{fig: outcomes} for some examples. 

\subsection{Computing the outcome}
The decision problem related to the Maker-Breaker domination game is the following.

\begin{problem}[\MBDG]\label{problem MB dom game}{}
\hspace{.2cm}

  \noindent  Input: a graph $G$, a first player Dominator or Staller. 
        
 \noindent   
 Question: what is the outcome of the Maker-Breaker domination game played on $G$ with the considered player going first? 
\end{problem}

\MBDG has been proved to be \PSPACE-complete even for restricted classes of graphs:

\begin{theorem}[\cite{makerbreaker, Oijid2025}]\label{thm:pspace}
    \MBDG is \PSPACE-complete, even restricted to bipartite graphs, split graphs or bounded degree graphs.
\end{theorem}

In this paper, we consider the parameterized complexity of this problem using several parameters. We refer the reader to standard books of parameterized complexity for definitions and general framework (see for example \cite{downey1999,Gurevich1984}). Note that hypothesis on the size of a dominating set of $G$ would not be enough to determine the outcome of an instance of the Maker-Breaker domination game. Indeed, even if $G$ has a dominating set of size 1, determining if the outcome is $\oN$ or $\oD$ is still hard. It cannot be $\oS$, since Dominator wins in one move when she starts.

\begin{proposition}\label{cor:dom2}
    \MBDG is \PSPACE-complete, even restricted to graphs $G$ with $\gamma(G)= 1$.
\end{proposition}

\begin{proof}
    We reduce our problem from \MBDG in a general graph. Let $G$ be a graph. 
 Let $G'$ be the graph obtained from $G$ by adding a universal vertex $v_0$ to $G$, thus $\gamma(G') = 1$. 
    Then Staller wins playing first in $G'$ if and only if Dominator wins playing first in $G$.
    Indeed, if Staller starts in $G'$, he must claim $v_0$ first, otherwise Dominator will claim it and win at her first move. Therefore, after this move, the game is played on $G'$ as it was played on $G$ with Dominator playing first. Note that $v_0$ will be dominated by the first move of Dominator. 
\end{proof}

\subsection{Union and join}

In~\cite{makerbreaker}, it is proved that \MBDG can be computed in linear time on cographs. This proof relies on the following lemma that computes the outcome of a disjoint union or a join of two graphs using the outcomes of both graphs.

\begin{lemma}[\protect{\cite{makerbreaker}}]\label{lemma union} 
    Let $G$ and $H$ be two graphs.
\begin{itemize}
\item If $o(G) = \oS$ or $o(H) = \oS$ then $o(G \cup H) = \oS$.
\item If $o(G) = o(H) = \oN$ then $o(G\cup H) = \oS$.
\item If $o(G) = o(H) = \oDom$ then $o(G \cup H) = \oDom$.
\item Otherwise, $o(G \cup H) = \oN$. 
\end{itemize}

\smallskip

\end{lemma}

\begin{lemma}[\protect{\cite{makerbreaker}}]\label{lemma join}

Let $G$ and $H$ be two graphs.

\begin{itemize}
\item If $G = K_1$ and $o(H) = \oS$ (or $H = K_1$ and $o(G) = \oS$), then $o(G \join H) = \oN$.
\item Otherwise, $o(G \join H) = \oDom$.
\end{itemize}
    
\end{lemma}

\subsection{Simplifying positions}

An important property of positional games is that giving a vertex to a player can only benefit him. This leads to the following observation:

\begin{observation}\label{obs: free vertex}
Let $P=(G, D, S)$ be a position and $v \in V(G)\setminus D\cup S$ be an unclaimed vertex of $P$. If Dominator wins on $(G,D,S)$, she wins on $(G, D \cup \{v\}, S)$. Similarly, if Staller wins on $(G,D,S)$, he wins on $(G, D, S \cup \{v\})$.
\end{observation}

To simplify some instances, one might need to remove vertices or edges of a position while keeping the same winning sets. This is formally expressed in the following observation.

\begin{observation}\label{obs:sameWS}
Let $P=(G, D, S)$ and $P'=(G', D', S')$ be two positions such that there exists a bijection $f:V(G)\setminus D\cup S \to V(G')\setminus D'\cup S'$ between the sets of unclaimed vertices. Assume that the vertices to complete a dominating set in $P$ and $P'$ are the same under the bijection $f$, i.e., $X\cup D$ is a dominating set of $G$ if and only if $f(X) \cup D'$ is a dominating set of $G'$. Then $P$ and $P'$ have the same outcome.
\end{observation}

Using Observation~\ref{obs:sameWS}, vertices claimed by Staller that are already dominated can be removed.


\begin{lemma}\label{removing vertices staller}
    Let $(G,D,S)$ be a position. Let $v \in S$ with a neighbor in $D$. Then $(G,D,S)$ and $(G\setminus\{v\},D,S\setminus\{v\})$ have the same outcome.
\end{lemma}

\begin{proof}
    Let $(G,D,S)$ be a position. Let $X \subset G\setminus(D\cup S)$ be a set of vertices. We have that $D \cup X$ is a dominating set of $G$ if and only if $X \cup D$ is a dominating set of $G\setminus\{v\}$. Therefore, the winning sets of $(G,D,S)$ and $(G\setminus\{v\},D,S\setminus\{v\})$ are the same, so, by Observation~\ref{obs:sameWS} they have the same outcome. 
\end{proof}

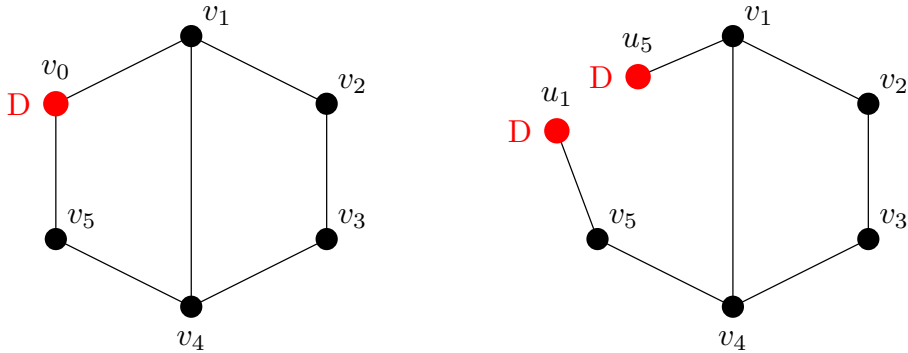
\begin{figure}[ht]
    \centering

\scalebox{1.2}{

\begin{tikzpicture}[scale=1.5]
  \coordinate (V1) at (0,1);
  \coordinate (V2) at (1,0.5);
  \coordinate (V3) at (1,-0.5);
  \coordinate (V4) at (0,-1);
  \coordinate (V5) at (-1,-0.5);
  \coordinate (V0) at (-1,0.5);

  \foreach \i/\label in {1/1, 2/2, 3/3, 5/5} {
    \draw (V\i) node[v] (1.5pt){} node[above right] {$v_{\label}$};
  }

  \draw (V1) -- (V2);
  \draw (V1) -- (V0);
  \draw (V2) -- (V3);
  \draw (V3) -- (V4);
  \draw (V4) -- (V5);
  \draw (V5) -- (V0);
  \draw (V1) -- (V4);

\draw (V0) node[R] (1.5pt){} node[above= .15cm] {$v_{0}$} node[left = .15cm] {\color{red} D};
\draw (V4) node[v] (1.5pt){} node[below = .15cm] {$v_{4}$};


  \coordinate (V1) at (4,1);
  \coordinate (V2) at (5,0.5);
  \coordinate (V3) at (5,-0.5);
  \coordinate (V4) at (4,-1);
  \coordinate (V5) at (3,-0.5);
  \coordinate (V01) at (2.7,0.3);
  \coordinate (V02) at (3.3,0.7);

  \foreach \i/\label in {1/1, 2/2, 3/3, 5/5} {
    \draw (V\i) node[v] (1.5pt){} node[above right] {$v_{\label}$};
  }

  \draw (V1) -- (V2);
  \draw (V1) -- (V02);
  \draw (V2) -- (V3);
  \draw (V3) -- (V4);
  \draw (V4) -- (V5);
  \draw (V5) -- (V01);
  \draw (V1) -- (V4);

\draw (V01) node[R] (1.5pt){} node[above= .15cm] {$u_{1}$} node[left = .15cm] {\color{red} D};
\draw (V02) node[R] (1.5pt){} node[above= .15cm] {$u_{5}$} node[left = .15cm] {\color{red} D};

\draw (V4) node[v] (1.5pt){} node[below= .15cm] {$v_{4}$};

\end{tikzpicture}
    }
    \caption{Split on the vertex $v_0$ using Lemma~\ref{removing vertices Dominator}. The two graphs have the same outcome.}
    \label{fig:split dom}
\end{figure}

The vertices claimed by Dominator cannot be removed in the same way, as they might be used to dominate other vertices that are currently unclaimed. However, the next lemma presents a way to transform instances where Dominator has claimed some vertices by {\em spliting} these vertices.  See Figure~\ref{fig:split dom} for an illustration.

\begin{lemma}\label{removing vertices Dominator}
    Let $(G,D,S)$ be a position. Let $v \in D$ of degree $d$, and let $v_1, \dots, v_d$ be the neighbors of $v$. Let $G'$ be the graph obtained from $G\setminus\{v\}$ and by adding $d$ new vertices $u_1, \dots, u_d$, private neighbors of $v_1, \dots, v_d$. Then $(G,D,S)$ and $(G',D\cup\{u_1,\dots,u_d\}\setminus\{v\},S)$ have the same outcome. If this happens, we say that we {\em split} the vertex $v$.
\end{lemma}

\begin{proof}
    Let $(G,D,S)$ be a position. Let $X \subset G\setminus(D\cup S)$ be a set of vertices. We have that $D \cup X$ is a dominating set of $G$ if and only if $X \cup G \cup \{u_1, \dots, u_d\} \setminus\{v\}$ is a dominating set of $G'$. Therefore, the winning sets of $(G,D,S)$ and $(G',D\cup\{u_1,\dots,u_d\}\setminus\{v\},S)$ are the same, so, by Observation~\ref{obs:sameWS} they have the same outcome. 
\end{proof}

\subsection{Some basic strategies}

 Leaves or support vertices play a special role in the Maker-Breaker domination game. They can be used to force moves. The following lemma shows that it is always in Staller's interest to play a support vertex, forcing Dominator to answer by claiming the leaf adjacent to it. 
 
\begin{lemma}[\protect{\cite[Proposition 20]{makerbreaker}}] \label{unclaimed tail}
    Let $(G,D,S)$ be a position. Let $v_0$ be an unclaimed leaf of $G$ and $v_1$ its support vertex. If it is Staller's turn and $v_1$ is unclaimed, then claiming $v_1$ is an optimal move. In particular, positions $(G,D,S)$ and $(G\setminus \{v_0,v_1\}, D, S)$ have the same outcome.
\end{lemma}

This lemma can be generalized when vertices have closed neighborhoods which are included one in the other. In this case, both players prefer to claim the vertex with the largest neighborhood.
This result is stated only for Dominator in~\cite{Intervalle}, but the same proof also works for Staller.

\begin{lemma}[\protect{\cite[Lemma 2.6]{Intervalle}}]\label{lemma: dominated move Maker-Breaker}
    Let $G=(V,E)$ be a graph, and $(G,D,S)$ be a position on $G$. Let $u,v \in V \setminus (D \cup S)$ be two unclaimed vertices of $G$ such that $N[u] \subset N[v]$. For any position of the game on which $u$ and $v$ are unclaimed, if a player has a winning strategy playing $u$, he has one playing $v$. 
\end{lemma}


A basic strategy that is often used is the {\em pairing strategy}. It has been used
to solve \MBDG in trees \cite{makerbreaker} and to prove that it is in NP for interval graphs \cite{Intervalle}. A pairing strategy consists of selecting disjoint pairs of vertices and, when the opponent claims one vertex of a pair, taking the other vertex of the pair. In this way, one can guarantee to take at least one vertex from each pair.


Finally, the next lemma, also called as the ``Super Lemma", states that if two unclaimed vertices behave in the same way in a given position, one can assume that each player will claim exactly one of the two vertices and that they can be claimed immediately. This will in particular be the case in the starting position if two vertices have the same closed neighborhood. This lemma has been stated in a general version for positional games in Oijid's PhD thesis~\cite{nacimthesis}, and is stated here for the Maker-Breaker domination game. 

\begin{lemma}[\protect{\cite[Super Lemma]{nacimthesis}}]\label{super lemma}
Let $(G,D,S)$ be a position of a Maker-Breaker domination game. Let $V' = V(G) \setminus (D\cup S)$ be the set of unclaimed vertices. Let $u,v \in V'$ be two vertices of $G$ such that, for any set $X \subset V'\setminus \{u,v\}$, we have $X \cup D \cup \{u\}$ is a dominating set of $G$ if and only if $X \cup D \cup \{v\}$ is a dominating set of $G$.

Then $(G,D,S)$ and $(G, D \cup \{u\}, S \cup \{v\})$ have the same outcome.    
\end{lemma}

\section{Number of moves}\label{sec:moves}

In this section, we consider as a parameter the number of moves a player needs to win. We will consider the more general context of Maker-Breaker positional games. Considering the number of moves as a parameter changes the nature of the problems we are considering. Since after $k$ moves,  it is possible that no player has won, we will consider the following problems.

\begin{problem}[\shortDD]\label{short-dom-dom}{}
\hspace{.2cm}

  \noindent  \Input: A graph $G$, the first player (Dominator or Staller),
  and an integer~$k$. 
    
 \noindent   \Parameter: $k$. 
    
 \noindent   
 \Question: Does Dominator have a strategy that claims a dominating set of $G$ in at most $k$ moves?
\end{problem}
\begin{problem}[\shortDS]\label{short-dom-stal}{}
\hspace{.2cm}

  \noindent  \Input: A graph $G$, the first player (Dominator or Staller),
  and an integer~$k$. 
    
  \noindent  \Parameter: $k$. 
    
 \noindent   
 \Question: Does Staller have a strategy that claims the neighborhood of a vertex in $G$ in at most $k$ moves?
\end{problem}


\begin{problem}[\shortMBM]\label{short-mb-maker}{}
\hspace{.2cm}

\noindent    \Input: A hypergraph $\hyp$, the first player (Maker or Breaker), 
and an integer~$k$. 
    
\noindent    \Parameter: $k$. 
    
\noindent    
\Question: Does Maker have a strategy to fill up a hyperedge of $\hyp$ in at most $k$ moves?
\end{problem}
\begin{problem}[\shortMBB]\label{short-mb-breaker}{}
\hspace{.2cm}

\noindent \Input: A hypergraph $\hyp$, the first player (Maker or Breaker), 
and an integer~$k$. 
    
\noindent    \Parameter: $k$. 
    
\noindent   
\Question: Does Breaker have a strategy to claim a transversal of $\hyp$ in at most $k$ moves?
\end{problem}

Among these four problems, the only one that has been studied is {\shortMBM}, proved to be \W[1]-complete by Bonnet \etal~\cite{Bonnet2017}.

\begin{theorem}[\protect Bonnet \etal~\cite{Bonnet2017}] \label{thm: short Maker-Breaker Bonnet}
    \shortMBM is \W[1]-complete.
\end{theorem}

We give the complexity of the three other problems in the rest of this section. 


\subsection{\shortDS is \W[1]-complete}

Since $\hypneighbor(G)$, the hypergraph of the closed neighborhoods of $G$, has a polynomial size in the size of $G$, the fact that \shortMBM is in $\W[1]$ implies that {\shortDS} is in~\W[1]. We prove that it is even \W[1]-complete.

\begin{theorem}\label{thm short dom game by staller}
    {\shortDS} is \W[1]-complete.
\end{theorem}

\begin{proof}
For the membership in~\W[1], we consider the Maker-Breaker domination game as a Maker-Breaker positional game with Staller playing the role of Maker. Let $G$ be a graph. The hypergraph of the positional game is then the hypergraph of the closed neighborhoods, $\hypneighbor(G)$ which has $|V(G)|$ vertices and at most $|V(G)|$ hyperedges.
Thus, it has a polynomial number of vertices and hyperedges. Since \shortMBM is in \W[1], deciding if Staller can win in $k$ moves is also in \W[1].

We now prove the completeness with a reduction from \shortMBM using ideas similar to those used in \cite{makerbreaker} to prove that \MBDG is \PSPACE-complete.

Let $\hxf$ be a hypergraph and $k$ be an integer. 
    Without loss of generality, suppose that every vertex in $\hyp$ is in at least one hyperedge. We construct a graph $G_{\hyp}=(V,E)$ as follows.

\begin{itemize}
    \item For each vertex $u_i$ in $\som$, we add a vertex $v_i$ in $V$. 
    \item For each hyperedge $f$ in $\WS$, we add $k+2$ new vertices $v_f^1, \dots, v_f^{k+2}$ in $V$. 
    \item For each pair of vertices $u_i, u_j \in \som$, we add an edge between $v_i$ and $v_j$ in $E$. 
    \item If a vertex $u_i$ of $\som$ belongs to a hyperedge $f$ of $\WS$, we add the edges $v_iv_f^1), \dots, v_iv_f^{k+2}$ to $E$.
\end{itemize}

This reduction is depicted in Figure~\ref{fig: reduction W1 MBDOM}.

\begin{figure}[t!]
    \centering
    \begin{subfigure}[t]{0.49\textwidth}
        \centering

\begin{tikzpicture}

\draw (0,0) node[v] (a) {} node[below = .2]{$u_3$};
\draw (1,0) node[v] (b) {} node[below = .2]{$u_4$};
\draw (0,1) node[v] (c) {} node[above = .2]{$u_1$};
\draw (1,1) node[v] (d) {} node[above = .2]{$u_2$};

\draw \hedgeii{a}{b}{2mm} node[below left = .5]{$e$};
\draw \hedgeiii{a}{c}{d}{2.5mm} node[below left = .2]{$f$};

\end{tikzpicture}
        
        \caption{A hypergraph $\hyp$}
    \end{subfigure}%
    ~ 
    \begin{subfigure}[t]{0.49\textwidth}
        \centering

\begin{tikzpicture}

\draw (-.375,2.4) node[v] (a) {} node[above = .2]{$v_1$};
\draw (1.125,1) node[v] (b) {} node[above = .15]{$v_2$};
\draw (2.625,1) node[v] (c) {} node[above = .15]{$v_3$};
\draw (4.125,2.4) node[v] (d) {} node[above = .2]{$v_4$};

\draw (a) -- (b);
\draw (a) -- (c);
\draw (a) -- (d);
\draw (b) -- (c);
\draw (b) -- (d);
\draw (c) -- (d);

\draw (-.75,0) node[v] (e11) {} node[below = .2]{$v_f^1$} ;
\draw (0,0) node[v] (e12) {} node[below = .2]{$v_f^2$} ;
\draw (.75,0) node[v] (e13) {} node[below = .2]{$v_f^3$} ;
\draw (1.5,0) node[v] (e14) {} node[below = .2]{$v_f^4$} ;

\draw (2.25,0) node[v] (e21)  {} node[below = .2]{$v_e^1$} ;
\draw (3,0) node[v] (e22) {} node[below = .2]{$v_e^2$} ;

\draw (3.75,0) node[v] (e23)  {} node[below = .2]{$v_e^3$} ;
\draw (4.5,0) node[v] (e24) {} node[below = .2]{$v_e^4$} ;

\draw (a) -- (e11);
\draw (b) -- (e11);
\draw (c) -- (e11);
\draw (a) -- (e12);
\draw (b) -- (e12);
\draw (c) -- (e12);
\draw (a) -- (e13);
\draw (b) -- (e13);
\draw (c) -- (e13);
\draw (a) -- (e14);
\draw (b) -- (e14);
\draw (c) -- (e14);

\draw (c) -- (e21);
\draw (d) -- (e21);
\draw (c) -- (e22);
\draw (d) -- (e22);
\draw (c) -- (e23);
\draw (d) -- (e23);
\draw (c) -- (e24);
\draw (d) -- (e24);

\end{tikzpicture}
        \centering

        \caption{The resulting graph.}
    \end{subfigure}
    \caption{The reduction in Theorem~\ref{thm short dom game by staller}, with $k= 2$}
    \label{fig: reduction W1 MBDOM}
\end{figure}
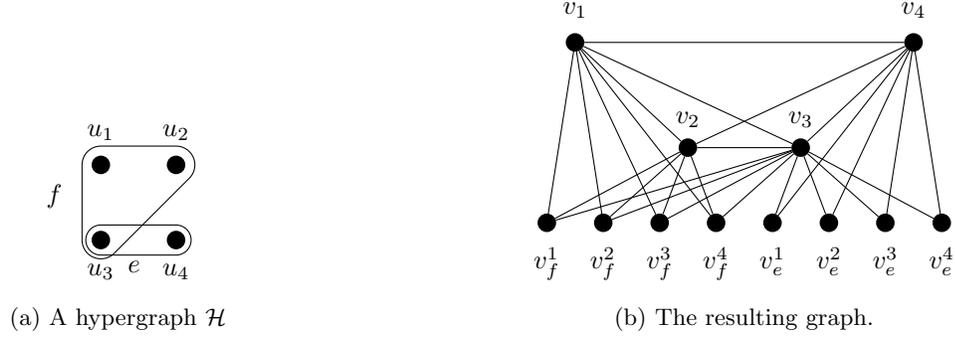

We prove that Staller can isolate a vertex of $G_{\hyp}$ in $k+1$ moves playing first (respectively second) if and only if Maker wins on $\hyp$ in $k$ moves playing first (respectively second).

Suppose that Maker has a winning strategy $\strat$ on $\hyp$ in $k$ moves playing first (respecively second). We provide a winning strategy for Staller in $k+1$ moves playing first (resp. second) on $G_{\hyp}$ as follows: 
\begin{itemize}
    \item If Maker goes first in $\hyp$ and claims $u_i$, Staller claims first the vertex $v_i$.
    \item If Dominator claims a vertex $v_i$, Staller answers by claiming the vertex $v_j$ such that $u_j$ is the answer of Maker when Breaker has claimed the vertex $u_i$ in $\strat$.
    \item If Dominator claims a vertex $v_f^i$ for $i \in \{1, ...,k+2\}$, Staller considers that Breaker has claimed an arbitrary unclaimed vertex $u_j$, and claims a vertex $v_j$ according to this move in $\strat$.
    \item When Staller have claimed $k$ vertices, as $\strat$ is a winning strategy in $\hyp$, there exists a hyperedge $f \in \WS$ fully claimed by Maker. Staller claims with his $(k+1)^{th}$ move a vertex $v_f^i$ with $1 \le i \le k+2$. Note that Dominator has claimed at most $k+1$ of them ($k$ if Staller was first to play on $\hyp$), thus one of them is available.
\end{itemize}

By construction, Staller with this strategy has claimed one vertex $v_f^i$ and all its neighbors, thus has won in $k+1$ moves.

\medskip

Suppose now that Breaker has a strategy $\strat$ in $\hyp$ preventing Maker from filling up a hyperedge in $k$ moves. We provide a strategy for Dominator in $G_{\hyp}$ preventin Staller to isolate a vertex in $k+1$ moves as follows:
\begin{itemize}
    \item If Breaker goes first in $\hyp$, Dominator claims first the vertex $v_i$ such that $u_i$ is the first vertex claimed in $\strat$.
    \item If Staller claims a vertex $v_i$, Breaker answers by claiming the vertex $v_j$ such that $u_j$ is Dominator's answer to the vertex $u_i$  claimed by Staller in $\strat$.
    \item If Staller claims a vertex $v_f^i$ for $i \in \{1, ...,k+2\}$, Dominator claims an arbitrary vertex of the graph, and considers that Maker has not played this move yet.
    \item When Staller claims his $k^{th}$ vertex, if she has claimed no vertex $v_f^i$, Dominator claim any unclaimed vertex of the graph. If she has claimed at least one vertex $v_f^i$, she claims the vertex $v_j$ corresponding to the vertex $u_j$ that would have been claimed according to $\strat$, which exists as the move on $v_f^i$ was ignored in $\strat$.
\end{itemize}

Following this strategy, if Staller did not play a vertex $v_f^i$ during his $k$ first moves, he cannot have claimed all the vertices corresponding to a hyperedge of $\hyp$. Therefore, he cannot isolate a vertex with his $(k+1)^{th}$ move. Since the vertices $v_i$ have at least $k+2$ neighbors each (we assumed that each vertex is in at least one hyperedge), and even if his last move would fill up a hyperedge, he has not claimed any vertex $v_f^i$, and therefore cannot have isolated any of them. If Staller has claimed at least one vertex $v_f^i$, she has claimed at most $k$ vertices $v_j$. Since Dominator's answers followed $\strat$, Staller cannot have claimed all the neighbors of a vertex $v_f^i$. In both cases, no vertex was isolated in the first $k+1$ moves of Staller, and Dominator has won. This reduction is clearly polynomial.
\end{proof}

\subsection{\shortMBB is in \W[2]}

Considering the number of moves required for Dominator to win is more complicated. Contrary to $\hypneighbor(G)$, the hypergraph of the dominating sets $\hypdom(G)$ does not have polynomial size in general. In particular, the known results for Maker-Breaker positional games when limiting the number of moves of Maker cannot be applied. We need to consider limited number of moves for Breaker. 

We prove in this subsection that the problem \shortMBB is in~\W[2]. The proof relies on an idea provided by Bonnet \etal~\cite{Bonnet2017}. In short, we first introduce 
a class of formulas $\MBWun$
in which all the $\forall$-quantifiers have their variable appearing only in inequalities except the last one, similar to the class $\forall^{\neq}$-$\FO$ from~\cite{Bonnet2017}, in which all the $\forall$-quantifiers have their variable appearing only in inequalities.
Then we provide an \FPT\ reduction from formulas in $\MBWun$ to formulas in $\Sigma_{2,1}$.  
We now brieﬂy review some standard definitions from mathematical logic.

 A {\em vocabulary} $\tau$ is a set of relational symbols
, each of which has a specified
 arity. 
A {\em $\tau$-structure} 
$  \alphabet = (A,	( R^{\alphabet} )_{R \in \tau})$
consists of a set $A$ together with an interpretation of
 each $k$-ary relational symbol $R$ from $\tau$ as a $k$-ary relation on $A$; that is, a
set $P^{\alphabet} \subseteq A^k$.
Atomic formulas over a vocabulary $\tau$ are of the form $x_1 = x_2$ or $R(x_1 , \ldots , x_k )$  where $R \in \tau $  and $x_1 , \ldots , x_k$ are variables. 
%


The class $\Sigma_1$ is the class of first-order formulas that can be written using only $\exists$-quantifiers and all at the beginning of the formula, i.e. formulas of the form $\exists x_1 \ldots  \exists x_k ~ \psi$ where $\psi$ is a quantifier-free formula.
The class $\Sigma_2$ is the class of first-order formulas that can be written with one block of consecutive $\exists$-quantifiers 
followed by one block of  $\forall$-quantifiers, i.e. formulas of the form $\exists x_1 \ldots  \exists x_k\forall y_1\ldots \forall y_{\ell} ~ \psi$ where $\psi$ is a quantifier-free formula.
In addition, for all $u\geq 1$, we denote by $\Sigma_{2,u}$ the class of $\Sigma_2$ 
formulas such that the last block of $\forall$-quantifiers has length at most~$u$.
For a given class of formulas~$\Phi$, the parameterized variant of model checking can be expressed as the following problem.

\begin{problem}[$MC(\Phi)$]\label{model-checking}{}
\hspace{.2cm}

\noindent \Input: A finite structure $\alphabet$ and a
formula $\phi \in \Phi$.

\noindent    \Parameter: $|\phi|$. 
    
\noindent   
\Question: 
Is $\phi$ satisfiable ?
\end{problem}
%
%

This parameterized problem can be useful to reason about parameterized complexity classes, since, for instance, $MC(\Sigma_1)$ is \W[1]-complete and $MC(\FO)$ is \AW[*]-complete~\cite{flum2003}. 
Furthermore, in~\cite{Bonnet2017}, Bonnet \etal introduce a new class of formulas by enriching $\Sigma_1$, allowing the addition of any number of $\forall$-quantifiers whose corresponding universal variables only appear in inequalities.
This new class, denoted by $\forall^{\neq}$-$\Sigma_1$, is then used to 
prove that {\shortMBM} is \W[1]-complete by a reduction to $MC(\Sigma_1)$.
Formally, $\forall^{\neq}$-$\Sigma_1$ is the
class of all first-order formulas of the form 
\[ Q_1 x_1 \dots Q_n x_n~\psi \] with
$Q_i\in \{\exists,\forall\}$,  $\psi$ a formula, every $\forall$-quantified variable $x_i$ only occurs in inequalities: $x_i \neq x_{j}$ that cannot be negated ?.

We extend $\Sigma_2$ in a similar way in order to prove 
that \shortMBB is in \W[2]. Let $\MBW$ be the class of formulas of the form \[ Q_1 x_1 \dots Q_n x_n \forall y_1 \dots \forall y_m~\psi \] where $\psi$ is a quantifier-free formula, $Q_i\in \{\exists,\forall\}$ and every $\forall$-quantified variable $x_i$ only occurs in inequalities: $x_i \neq x_{j}$, and these inequalities cannot be negated.  Note that we \textbf{do not} allow inequalities of the form $x_i \neq y_j$ where $y_j$ is quantified with a $\forall$.
Finally, let $\MBWU$ the class of $\Sigma_{2,u}$ where the last block of $\forall$-quantifiers has length at most $u$.

Bonnet \etal~\cite{Bonnet2017} prove that $MC(\forall^{\neq}$-$\Sigma_1)$
is in \W[1] by reducing formulas of the new class to formulas in $\Sigma_1$ using an {\FPT} reduction. We will similarly reduce formulas in $\MBWU$ 
to formulas in $\Sigma_{2,u}$ to prove that $MC(\MBWU)$ is in \W[2].
We then use a result from Downey \etal~\cite{Downey1996} and Flum and Grohe~\cite{flum2003}, stating that 
for any $u\geq 1$,
model checking $\Sigma_{2,u}$, parameterized by the size of the formula, is \W[2]-complete. 

\begin{lemma} \label{model checking W2}
   $\forall u \geq 1$,  $MC(\MBWU)$ is \W[2]-complete parameterized by the size of the formula.
\end{lemma}


\begin{proof}
The hardness follows from the $\W[2]$-hardness of model checking  $\Sigma_{2,u}$, as $\Sigma_{2,u} \subset \MBWU$.
For the membership, let us fix some $u\geq 1$ and let $(\alphabet, \phi)$ be an instance of $MC(\MBWU)$. 
We reduce it to an instance of $MC(\Sigma_{2,u})$.
If $\phi$ contains no $\forall$-quantifier followed by an $\exists$-quantifier,
then $(\alphabet, \phi)$ is already an instance of $MC(\Sigma_{2,u})$. 
Hence, we can write $\phi$ as: 
\[\phi = Q_1 x_1 Q_2 x_2 \ldots Q_{i-1}x_{i-1} \forall x_i \exists x_{i+1}
\ldots \exists x_{k} \forall y_1 \dots \forall y_{\ell}~ \psi\] with $x_i$ the rightmost universal quantified variable that is followed by an 
existential variable and with $\ell \leq u$.

The proof follows the same pattern as the one by Bonnet \etal~\cite{Bonnet2017}, only, instead of eliminating all universal quantifiers until ending up with a $\Sigma_1$ formula, we 
will eliminate all the universal quantifiers 
up until $x_i$'s, ending up with a $\Sigma_{2,u}$ formula. 
We will first show how to eliminate the universal quantifier of $x_i$.
%
%
The idea is that, since $x_i$ only appears in inequalities with other variables $x_j$s, we only need to distinguish whether the value assigned to $x_i$ is equal to that assigned to some other variable or not. Indeed, if it is not, i.e. if the formula is satisfied for an assignment of $x_i$ such that all inequalities between $x_i$ and other variables are satisfied, then it is 
still satisfied if $x_i$'s assignment is changed but all these inequalities still hold.
The technique will then be used iteratively to eliminate all the universal quantifiers before $x_i$,
recall that the corresponding variables only appear in inequalities.

We will prove that we can replace
the subformula 

\noindent $\phi_1(x_1, \ldots, x_{i-1}) = \forall x_i \exists x_{i+1} \ldots \exists x_k \forall y_1 \dots \forall y_{\ell} ~  \psi $ by

$$\phi_2(x_1, \ldots, x_{i-1}) = $$

\vspace{-0.9cm}
\begin{eqnarray}
 &\exists z_i \exists z_{i+1} \ldots \exists z_{k}
\big( ( \forall y_1 \dots \forall y_{\ell} ~ \psi[ z_i/x_i, z_{i+1}/x_{i+1}, \ldots, z_k/x_k
] ) \wedge \phantom{phantomminmax} \\
& \bigwedge_{j=1}^{i-1}  \exists z_{i+1}^j \ldots \exists z_{k}^j \forall y_1 \dots \forall y_{\ell} ~ \psi[ x_j/x_i, z_{i+1}^j/x_{i+1}, \ldots, z_k^j/x_k 
 ] \wedge \\
&  \bigwedge_{j=i+1}^{k}  \exists z_{i+1}^j \ldots \exists z_{k}^j \forall y_1 \dots \forall y_{\ell} ~ \psi[ z_j/x_i, z_{i+1}^j/x_{i+1}, \ldots, z_k^j/x_k 
  ]  \big)
\end{eqnarray}

where $\phi[y/x]$  means that in $\phi$ every occurrence of $x$ is replaced by $y$.
This reduction is an \FPT-reduction since the size of the formula $\phi_2$ is a function of the size of the formula $\phi_1$. We prove that, for $c_1, \dots, c_{i-1}$ elements of $\alphabet$, $\phi_1(c_1, \dots, c_{i-1})$ is True if and only if $\phi_2(c_1, \dots, c_{i-1})$ is.

We first give the intuition about the construction of $\phi_2$. The $\forall x_i$ is split depending on the value assigned to $x_i$: 
\begin{itemize}
    \item $(1)$ is here to ensure that at least one valuation satisfies $\psi$.
    \item $(2)$ is here to verify that for all the values of $x_j$, with $j < i$ that are already chosen, $\psi$ can still be satisfied if $x_i$ takes one of these values.
    \item $(3)$ ensures that, whenever $x_i$ takes a value that was required for our first valuation that satisfies $\psi$ (in $(1)$), it is possible to find another valuation for the $x_j$s with $j > i$ that will satisfy $\psi$.
\end{itemize}
Note that, since $x_i$ only appears in inequalities, these are the only cases where a value of $x_i$ can turn the value of $\psi$ to False.

First if  $\phi_1(c_1, \ldots, c_{i-1})$ is True then  $\phi_2(c_1, \ldots, c_{i-1})$ is, since $(1), (2)$  and $(3)$ are special cases of the universal quantification over $x_i$. Indeed, if  $\phi_1(c_1, \ldots, c_{i-1})$ is True, given any assignment $c_i \in A$ of $x_i$, there exists an assignment $ C = c_{i+1} \ldots c_{k}$ of $x_{i+1} \ldots x_{k}$ such that $\psi$ is true (and thus $(1)$ is True). In particular, this is also the case for assignments $c_i$ equals to assignments $c_j$ with $j < i$ (meaning $(2)$ is True) and this is also the case for assignments $c_i$ equals to any $c_j$ defined in $C$ with $j > i$ i.e. there still exists in this case an assignment $c'_{i+1} \ldots c'_{k}$ of $x_{i+1} \ldots x_{k}$ such that $\psi$ is true (meaning $(3)$ is True).


Suppose now that $\phi_1(c_1, \dots, c_{i-1})$ is False. Suppose by contradiction that $\phi_2(c_1, \dots, c_{i-1})$ is True, i.e. there is an assignment $c_i, c_{i+1}, \ldots c_{k}$ of $z_i, z_{i+1}, \ldots z_{k}$ such that the remaining subformula is True.
By hypothesis, there exists $c_i' \in \alphabet$ such that $\phi_1(c_1, \dots, c_{i-1})[c_i'/x_i]$ is False. 
\begin{itemize}
    \item If $c_i' = c_j$ for some $1 \le j \le i-1$. Let  $c_{i+1}', \dots, c_k' \in \alphabet$ be assignments for the variables $z_{i+1}^j, \dots, z_k^j$ which satisfy $\phi_2$, i.e. such the subformula corresponding to the $j^{th}$ term in $(2)$, $\forall y_1 \dots \forall y_{\ell}~ \psi[x_j/x_i, c_{i+1}/x_{i+1}, \dots c_k/x_k]$, is True. Assigning $c_{j}, c_{i+1}', \dots, c_k'$ to  the variables $x_{i}, x_{i+1}, \dots, x_k$ in $\phi_1$ would thus make the formula $\phi_1(c_1, \dots, c_{i-1})$ hold when assigning $c_i$ to $x_i$, which contradicts the hypothesis.
    \item If $c_i' = c_j$ for some $i+1 \le  j \le k$. Let $c_{i+1}', \dots, c_k'$ be assignments for the variables $z_{i+1}^j, \dots z_k^j$  which satisfy $\phi_2$, i.e. such the subformula corresponding to the $j^{th}$ term in $(3)$,
    $\forall y_1 \dots \forall y_{\ell}~ \psi[z_j/x_i, c_{i+1}/x_{i+1}, \dots c_k/x_k]$, is True. Assigning $c_{j}, c_{i+1}', \dots, c_k'$ to  the variables $x_{i}, x_{i+1}, \dots, x_k$ in $\phi_1$ would thus make the formula $\phi_1(c_1, \dots, c_{i-1})$ hold when assigning $c_i$ to $x_i$, which contradicts the hypothesis.
    \item Otherwise, let $l_1, \dots, l_m$ be the literals that contain $x_i$ in $\psi$. Recall that all of them are inequalities of the form $x_i \neq x_j$ for $j \neq i$ . Recall $c_i$ is an assignment of $z_i$ which satisfies $\phi_2$, and let $l_1', \dots, l_m'$ be the literals of $\psi[c_i/x_i]$ for some valuation that satisfies $\phi_2$. 
    Let $l_1', \dots, l_m'$ be the corresponding literals 
    in $\forall y_1 \dots \forall y_{\ell} ~ \psi[ z_i/x_i, z_{i+1}/x_{i+1}, \ldots, z_k/x_k
]$. Remark that $\forall y_1 \dots \forall y_{\ell} ~ \psi[c_i/x_i, c_{i+1}/x_{i+1}, \dots, c_k/x_k]$ is True.
    Some of $l_1', \dots, l_m'$ might be False, but, as they cannot appear negated, turning them to True cannot turn the formula from True to False. Since changing the assignment of $x_i$ from $c_i$ to $c_i'$ renders all these inequalities True,
    $\forall y_1 \dots \forall y_{\ell} ~  \psi[c_i'/x_i, c_{i+1}/x_{i+1}, \dots, c_k/x_k]$ is True which again leads to a contradiction. 
\end{itemize}
This concludes the proof that the formulas are equivalent. The formula $\phi_2$ is not exactly in  $\Sigma_{2,u}$, since there are existential and universal quantifiers following a conjunction. However the formula can be transformed by moving all quantifiers right before the conjunctions. While there is an increase in the number of existential quantifiers compared to $\phi_1$, the universal quantifiers  $\forall y_1 \dots \forall y_{\ell}$ can be shared by all of the conjunctions, leading to a formula of the form

\begin{eqnarray*}
  \phi_2'(x_1, \ldots, x_{i-1}) = &\exists z_i \exists z_{i+1} \ldots \exists z_{k} \exists z_{i+1}^1 \ldots \exists z_{i+1}^k \ldots \exists z_{k}^k ~ \forall y_1 \dots \forall y_{\ell} \\
& \big( ( \psi[ z_i/x_i, z_{i+1}/x_{i+1}, \ldots, z_k/x_k
] ) \wedge \phantom{max} \\
& \bigwedge_{j=1}^{i-1} \psi[ x_j/x_i, z_{i+1}^j/x_{i+1}, \ldots, z_k^j/x_k, 
 ] \wedge \\
&  \bigwedge_{j=i+1}^{k} \psi[ z_j/x_i, z_{i+1}^j/x_{i+1}, \ldots, z_k^j/x_k, 
  ]  \big)
\end{eqnarray*}

belonging to $\MBWU$.
By iterating this construction until all the universal quantifiers before the last alternation are removed, we conclude that $MC(\MBWU)$ is in \W[2]. 
\end{proof}


\begin{corollary}\label{corollary MB param break in W2}
    \shortMBB\ is in \W[2].
\end{corollary}

\begin{proof}
    The winning condition for Breaker in \shortMBB can be expressed
    in $\MBWun$ as follows:
\[ \forall v_1 \in V ~ \! \exists u_1 \in V ~ \! \forall v_2 \in V ~ \! \dots \forall v_k \in V ~ \! \exists u_k \in V ~ \! \forall e \in E ~ \left ( \bigwedge_{i=1}^{k}  \bigwedge_{j=1}^{i} u_i \neq v_j \right )  \wedge \left ( \bigvee_{i=1}^{k}   Inc(u_i,e)  \right )  \]

where $Inc(u_i,e)$ is True if and only if $e$ is incident to $u_i$.

Therefore, by Lemma~\ref{model checking W2}, it is in \W[2] as the $v_i$'s only appear in inequalities.
\end{proof}

\subsection{\W[2]-completeness}

In the previous subsection we proved  that 
{\shortMBB}
is in \W[2]. As a consequence,
{\shortDD}
is in \W[2]. 
We prove that it is even \W[2]-complete. 


\begin{theorem}\label{shortDD W2 hard}
{\shortDD} is \W[2]-complete.
\end{theorem}

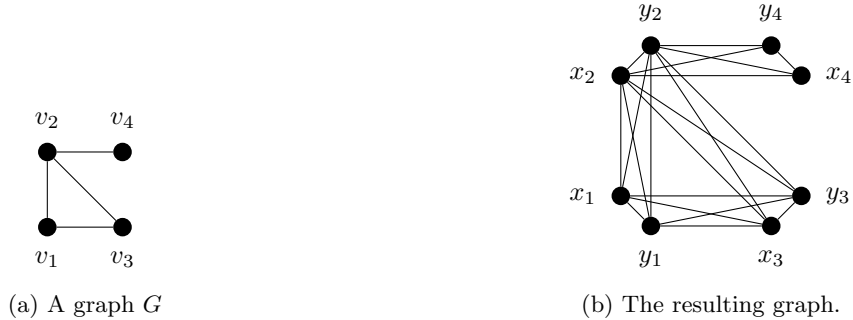
\begin{figure}
    \centering
    \begin{subfigure}[t]{0.49\textwidth}
        \centering

\begin{tikzpicture}

\draw (0,0) node[v] (a) {} node[below = .2]{$v_1$};
\draw (0,1) node[v] (b) {} node[above = .2]{$v_2$};
\draw (1,0) node[v] (c) {} node[below = .2]{$v_3$};
\draw (1,1) node[v] (d) {} node[above = .2]{$v_4$};

\draw (a) -- (b);
\draw (c) -- (b);
\draw (a) -- (c);
\draw (b) -- (d);

\end{tikzpicture}
        
        \caption{A graph $G$}
    \end{subfigure}%
    ~ 
    \begin{subfigure}[t]{0.49\textwidth}
        \centering

\begin{tikzpicture}
\draw (-0.2,0.2) node[v] (a1) {} node[left = .2]{$x_1$};
\draw (-0.2,1.8) node[v] (b1) {} node[left = .2]{$x_2$};
\draw (1.8,-0.2) node[v] (c1) {} node[below = .2]{$x_3$};
\draw (2.2,1.8) node[v] (d1) {} node[right = .2]{$x_4$};
\draw (0.2,-0.2) node[v] (a2) {} node[below = .2]{$y_1$};
\draw (0.2,2.2) node[v] (b2) {} node[above = .2]{$y_2$};
\draw (2.2,0.2) node[v] (c2) {} node[right = .2]{$y_3$};
\draw (1.8,2.2) node[v] (d2) {} node[above = .2]{$y_4$};

\draw (a1) -- (a2);
\draw (c2) -- (c1);
\draw (b2) -- (b1);
\draw (d2) -- (d1);

\draw (a1) -- (b1);
\draw (a1) -- (b2);
\draw (b1) -- (c1);
\draw (b1) -- (c2);
\draw (c1) -- (a1);
\draw (c1) -- (a2);

\draw (a2) -- (b1);
\draw (a2) -- (b2);
\draw (b2) -- (c1);
\draw (b2) -- (c2);
\draw (c2) -- (a1);
\draw (c2) -- (a2);

\draw (b1) -- (d1);
\draw (b2) -- (d1);
\draw (b1) -- (d2);
\draw (b2) -- (d2);

\end{tikzpicture}
        \centering

        \caption{The resulting graph.}
    \end{subfigure}
    \caption{The reduction in Theorem~\ref{shortDD W2 hard}}
    \label{fig: reduction W2hard MBDOM}
\end{figure}

\begin{proof}
To prove that {\shortDD} is in \W[2], we consider the Maker-Breaker domination game as a Maker-Breaker positional game with Dominator playing the role of Breaker and thus the closed neighborhood hypergraph. Since it has a polynomial number of vertices and hyperedges and since \shortMBB is in \W[2], deciding if Dominator can win in $k$ moves is also in \W[2].

We now prove that \shortDD is \W[2]-hard using a reduction from the following {\kdomset} problem
,
 which is known to be \W[2]-complete from Downey and Fellows~\cite{Downey1995}.
 
\begin{problem}[\kdomset]\label{kdomset}{}
\hspace{.2cm}

\noindent \Input: A graph $G=(V,E)$, an integer $k$.

\noindent    \Parameter: $k$. 
    
\noindent   \Question: Does there exist a dominating set $D\subseteq V$ of size $k$? 
\end{problem}


Let $(G, k)$ be an instance of the 
{\kdomset} problem. We consider the graph $G'$ with:
\begin{itemize}
    \item $V(G')= \{x_i \mid v_i \in V(G)\} \cup \{y_i \mid v_i \in V(G)\}$
    \item $E(G') = \{x_iy_i \mid v_i \in V(G)\} \cup \{x_iy_j \mid v_iv_j \in E(G) \} \cup  \{x_ix_j \mid v_iv_j \in E \} \cup \{y_iy_j \mid v_iv_j \in E(G) \} $.
\end{itemize}

The reduction is depicted in Figure~\ref{fig: reduction W2hard MBDOM}.

Intuitively, the graph $G'$ is obtained from $G$ by duplicating all the vertices of $G$. Therefore, if $G$ has a dominating set $D$ of size $k$, Dominator can construct a dominating set on $G'$ in $k$ moves using a pairing strategy with pairs $\{(x_i,y_i)\mid v_i \in D\}$ since any set containing one vertex from each pair is a dominating set.

\end{proof}

As a consequence, since {\shortDD} can be reduced in polynomial time to {\shortMBB}, 
{\shortMBB}
is also \W[2]-hard.


\begin{corollary}\label{shortmbb w2 complete}
 {\shortMBB} is \W[2]-complete.
\end{corollary}

The Corollary~\ref{shortmbb w2 complete} is quite surprising, since \shortMBM is \W[1]-complete. Moreover, this shows that even if the roles of Maker and Breaker can be exchanged by considering the transversal hypergraph, this has a cost in terms of complexity. This can be explained by the fact that the winning condition for Maker is local (i.e. filling up a hyperedge), while the winning condition for Breaker is global.

%
%
%

\

 On the positive side, since both verifying if there exists a dominating set or a closed neighborhood of a vertex fully claimed by a player can easily be expressed in first order logic using $O(k)$ quantifiers, meta-theorems on model checking from Grohe \etal~\cite{Grohe2013}, and Dreier \etal~\cite{Dreier2024} provides us the following theorem.

\begin{theorem}
{\shortDD} and {\shortDS} can be solved in \FPT\ time, parameterized by the number of moves on monadically stable graph classes.
\end{theorem}

Note that, since nowhere dense graphs and bounded treewidth graphs are monadically stable, \MBDG\ can be computed in \FPT\ time parameterized by the number of moves in these classes of graph.


\section{Graphs with modules}\label{sec:modules}

In the rest of the paper, we only consider the Maker-Breaker domination game with unlimited number of moves.
We will handle the parameterized complexity of the Maker-Breaker domination game considering some graph parameters. This section focuses on parameters related to modules. Being able to handle modules makes it possible to generalize the polynomial algorithm for cographs presented in~\cite{makerbreaker} to parameters which can be seen as some distance to cographs. In particular, Lemma~\ref{super lemma module edition} also provides \FPT\ algorithms for the modular-width and the $P_4$-fewness.

\subsection{Replacing a module in a graph}

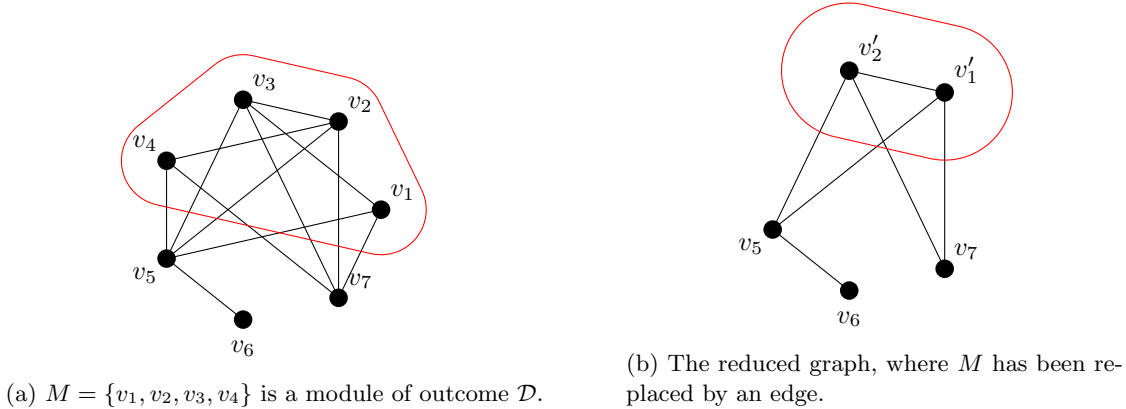
\begin{figure}[ht]
    \centering

\begin{subfigure}{.47 \textwidth}

\centering

    \begin{tikzpicture}[scale=1.5]
    \foreach \i/\label in {1/1, 2/2, 3/3, 4/4, 5/5, 6/6, 7/7} {
    \coordinate (V\i) at ({360/7 * (\i - 1)}:1);
    }
  \foreach \i/\label in {1/1, 2/2, 3/3, 7/7} {
    \draw (V\i) node[v] (1.5pt){} node[above right] {$v_{\label}$};
  }
    \draw (V5) node[v] (1.5pt){} node[below left] {$v_{5}$};
    \draw (V4) node[v] (1.5pt){} node[above left] {$v_{4}$};
    \draw (V6) node[v] (1.5pt){} node[below = .15cm] {$v_{6}$};

  \foreach \i in {1, 2, 3, 4} {
      \draw (V\i) -- (V7);
      \draw (V\i) -- (V5);
    }
  \draw (V5) -- (V6);
  \draw (V1) -- (V3);
  \draw (V4) -- (V2);
  \draw (V3) -- (V2);

\draw[color = red] \hedgeiiii{V3}{V2}{V1}{V4}{4mm};

\end{tikzpicture}

\caption{$M = \{v_1, v_2, v_3, v_4\}$ is a module of outcome~$\oDom$.}
    
\end{subfigure} \hfil \begin{subfigure}{.4\textwidth}

\centering

\begin{tikzpicture}[scale=1.5]

    \foreach \i/\label in { 2/2, 3/3, 5/5, 6/6, 7/7} {
    \coordinate (V\i) at ({360/7 * (\i - 1)}:1);
    }

\draw (V2) node[v] (1.5pt){} node[above right] {$v'_{1}$};
\draw (V3) node[v] (1.5pt){} node[above right] {$v'_{2}$};
\draw (V7) node[v] (1.5pt){} node[above right] {$v_{7}$};
\draw (V5) node[v] (1.5pt){} node[below left] {$v_{5}$};

    \draw (V6) node[v] (1.5pt){} node[below = .15cm] {$v_{6}$};

  \foreach \i in { 2, 3} {
      \draw (V\i) -- (V7);
      \draw (V\i) -- (V5);
    }
  \draw (V5) -- (V6);
  \draw (V3) -- (V2);

\draw[color = red] \hedgeii{V3}{V2}{6mm};

\end{tikzpicture}

\caption{The reduced graph, where $M$ has been replaced by an edge.}
    
\end{subfigure}
    
    \caption{Reduction of module using Lemma~\ref{super lemma module edition}}
    \label{fig:module reduction}
\end{figure}

Most results of this section are obtained thanks to the following lemma, which is a generalization of Lemma~\ref{super lemma}. 
It states that any module of size at least 3 can be replaced by a module of size at most~3 without changing the outcome of the graph. 
This lemma is illustrated in Figure~\ref{fig:module reduction}.

\begin{lemma}\label{super lemma module edition}

Let $G$ be a graph, and $M$ be a module of $G$ such that $|M| \ge 3$. Denote by $P_i$ the path on $i$ vertices.

\begin{itemize}
    \item If $o(G[M]) = \oDom$, let $G'$ be the graph obtained from $G$ by replacing $M$ by a copy of $P_2$. 
    \item If $o(G[M]) = \oS$, let $G'$ be the graph obtained from $G$ by replacing $M$ by two independent vertices.
    \item If $o(G[M]) = \oN$, let $G'$ be the graph obtained from $G$ by replacing $M$ by a copy of $P_3$.
\end{itemize}

Then $G$ and $G'$ have the same outcome.
\end{lemma}

\begin{proof}
    Let $G = (V,E)$ be a graph and let $M$ be a module of $G$. Suppose, for instance that $o(G[M]) = \oD$, the other cases being similar. Let $G'$ be the graph obtained from $G$ by replacing $M$ by a $P_2$. Let $\{u,v\}$ be the vertices of the $P_2$ in $G'$.

    First, we apply Lemma~\ref{super lemma} in $G'$ to $\{u,v\}$ that have the same neighborhood in $G'\setminus \{u,v\}$. Thus, $G'$ and $(G',\{u\},\{v\})$ have the same outcome.
    
    Suppose that Dominator has a winning strategy $\strat$ in $(G',\{u\},\{v\})$. We consider the following strategy for Dominator: she plays in $G\setminus M$ following $\strat$ and in $G[M]$ following her winning strategy as a second player (if one of the strategies is over, she plays randomly).

Following this strategy, we prove that Dominator will always dominate all the vertices of $G$. Let $x$ be a vertex of $G$.
\begin{itemize}
    \item If $x \in (V \setminus M) \cap N[M]$. Since $|M| \ge 2$, she will always claim at least one vertex $x_0$ of $M$ that will dominate $x$.
    \item If $x \in V \setminus N[M]$, let $y$ be the vertex that dominates $x$ in $G'$. This vertex exists since $\strat$ is a winning strategy in $G'$.  Since, $x$ is not adjacent to $M$, $y \neq u$. Thus, her strategy in $G$ also claims $y$ and $x$ is dominated.
    \item If $x \in M$, since $o(G[M]) = \oDom$ and she has played according to her winning strategy in $G[M]$, she will dominate $x$.
\end{itemize}

The same proof works if $O(G[M]) = \oN$ by applying Lemma~\ref{super lemma} to the two leaves of $P_3$ or if $o(G[M]) = \oS$ by applying Lemma~\ref{super lemma} to the two isolated vertices.

Suppose now that Staller has a winning strategy $\strat$ in $(G',\{u\},\{v\})$. Consider the strategy where Staller plays according to $\strat$ in $G \setminus M$ when Dominator plays in $G\setminus M$ and plays arbitrary vertices in $M$ otherwise.

For the same reason as previously, if Staller manages to isolate a vertex $x$ of $G\setminus M$ in $G'$ according to $\strat$, he still isolates it in $G$. Indeed, since $x$ is isolated in the game played on $(G',\{u\},\{v\})$, $x$ is not adjacent to $u$ in $G'$ and thus it has no neighbor in $M$ in $G$. Thus, Staller has a winning strategy in $G$.

If $o(G[M]) = \oN$ or $o(G[M]) = \oS$, the same argument works by adding the case that if Staller isolates a vertex of $P_3$ or of the two isolated vertices that substitutes $M$, he has claimed all the neighbors of $M$, and therefore has isolated a vertex of $M$ playing according to his winning strategy in $G[M]$.

Finally, we have proved that $G$ and $G'$ have the same outcome.
\end{proof}

This lemma directly gives a linear kernel for \MBDG parameterized by the {\em neighborhood diversity}, i.e. the minimum integer $w$ such that $V$ can be partitioned into $w$ sets $V_1, \dots, V_w$ and each $V_i$ is a set of vertices of same type, i.e. a module which induces either a clique or a stable set.

\begin{theorem} \label{theorem neighborhood diversity}
    \MBDG has a linear kernel, parameterized by the neighborhood diversity of the graph.
\end{theorem}

\begin{proof}
Let $G = (V,E)$ be a graph of neighborhood diversity $w$. Let $V_1, \dots, V_w$ be a partition of $V$ such that for $1 \le i \le w$, $V_i$ is a set of vertices of same type - i.e. if $x,y \in V_i$, $N(x) = N(y)$ or $N[x] = N[y]$. Therefore, by Lemma~\ref{super lemma module edition}, if $V_i$ contains at least two vertices, it can be replaced by either an edge or two independent vertices without changing the outcome. After applying this operation $w$ times (which can be done in linear time), we obtain a graph on at most $2w$ vertices and the same outcome, which gives a linear kernel.
\end{proof}

\subsection{Modular-width}\label{sec:neighbours}

Let $G_0 = (V_0, E_0)$ be a graph, and denote its vertices $v_1, \dots v_n$. Let $G_1 = (V_1, E_1), \dots, G_n = (V_n, E_n)$ be graphs. The {\em substitution operation} with respect to $G_0$, denoted by $G_0(G_1, \dots, G_n)$, consists in replacing each vertex $v_i \in V_0$ by the graph $G_i$ and adding all the edges between $G_i$ and $G_j$ if $v_iv_j \in E_0$. $G_i$ is then a module of the resulting graph.

\begin{figure}[ht]
    \centering
\begin{tikzpicture}
    \draw[style = dashed] (-6,0) ellipse (1cm and 2cm);
    \draw[style = dashed] (-2,0) ellipse (1cm and 2cm);
    \draw[style = dashed] (2,0) ellipse (1cm and 2cm);
    \draw[style = dashed] (6,0) ellipse (1cm and 1cm);

    \draw[style = dashed] (-6,-.5) ellipse (.5cm and 1cm);

  \draw[style = dashed, line width=1.8pt] (-5,0) -- (-3,0);
  \draw[style = dashed, line width=1.8pt] (-1,0) -- (1,0);
  \draw[style = dashed, line width=1.8pt] (3,0) -- (5,0);

\draw (-6,-1) node[v] (v1){} node[below = .15cm] {$v_{0}$};
\draw (-6,0) node[v] (v0){} node[left = .15cm] {$v_{1}$};
\draw (-6,1) node[v] (v2){} node[above = .15cm] {$v_{2}$};

\draw (-2,-1) node[v] (v3){} node[below = .15cm] {$v_{3}$};
\draw (-2,1) node[v] (v4){} node[above = .15cm] {$v_{4}$};

\draw (2,-1) node[v] (v5){} node[below = .15cm] {$v_{5}$};
\draw (2,1) node[v] (v6){} node[above = .15cm] {$v_{6}$};

\draw (6,0) node[v] (v7){} node[below = .15cm] {$v_{7}$};

\draw (v0) -- (v3);
\draw (v0) -- (v4);
\draw (v1) -- (v3);
\draw (v1) -- (v4);
\draw (v2) -- (v3);
\draw (v2) -- (v4);
\draw (v0) -- (v1);

\draw (v4) -- (v3);
\draw (v5) -- (v6);

\draw (v3) -- (v5);
\draw (v3) -- (v6);
\draw (v4) -- (v5);
\draw (v4) -- (v6);

\draw (v5) -- (v7);
\draw (v6) -- (v7);

\end{tikzpicture}    
\caption{Example of a modular decomposition of width $4$.  Here each set of vertices in circled with dashed lines are modules.}
    \label{fig:modular decomposition}
\end{figure}
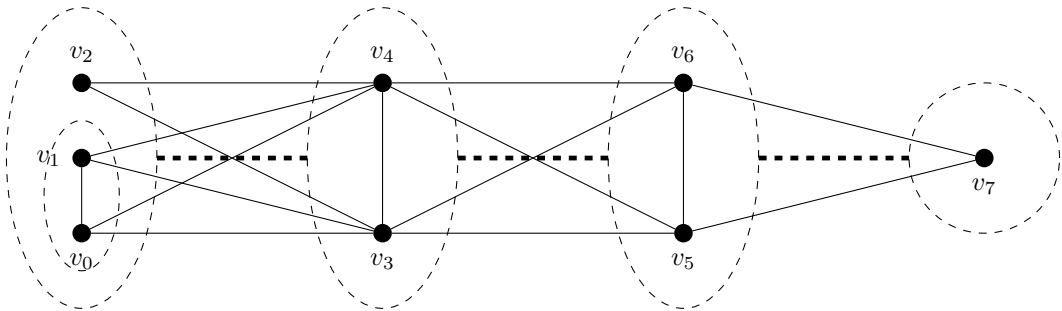

A {\em modular decomposition} of a graph $G = (V,E)$ is an expression that constructs $G$ using the following operations: addition of an isolated vertex, disjoint union of two graphs, complete join of two graphs, and substitution operation. 
    The {\em width} of a modular-decomposition is the {\em order} - i.e.  number of vertices - of the largest graph $G_0$ in the substitution operation. The {\em modular-width} of $G$, denoted by $mw(G)$ is the minimum integer $m$ such that $G$ can be obtained from a modular-decomposition of width $m$. An example of modular decomposition is provided in Figure~\ref{fig:modular decomposition}. 

\begin{theorem}
  \MBDG is \FPT{} parameterized by the modular-width of the graph.
\end{theorem}

\begin{proof}
    Let $G = (V,E)$ be a graph, and let $k = mw(G)$. First, we compute a decomposition of $G$ in, which can be done in linear time using an algorithm from Tedder \etal~\cite{Tedder2008}. We embed it in a tree such that, the nodes are the join, union and substitution operations, and the leaves are isolated vertices. Now we compute the outcome of $G$, with a bottom-up algorithm on the tree structure. Starting from the leaves, which all have outcome $\mathcal{N}$, we compute the outcome on each node, considering the subgraph $H$ induced by its descendant, as follows: 
\begin{enumerate}
    \item If $|V(H)| \le 3k$, we compute the winner by an exhaustive computation of all the possible games. 
    \item If $H = G_1 \cup G_2$, we compute inductively the outcome of $G_1$ and $G_2$ which have modular-width at most $k$ and we return the outcome according to Lemma~\ref{lemma union}.
    \item If $H = G_1 \join G_2$, we compute inductively the outcome of $G_1$ and $G_2$ which have modular-width at most $k$, and we return the outcome according to Lemma~\ref{lemma join}
    \item If $H = G'(G_1, \dots, G_p)$, necessarily, we have $p \le k$. We replace, according to Lemma~\ref{super lemma module edition}, each $G_i$ having more than three vertices by a $P_2$, a $P_3$ or two isolated vertices depending on its outcome. Then, we compute $o(H)$ as it has at most $3k$ vertices.
\end{enumerate}

This algorithm runs in time $O\left ( 3^{3k} n \right )$, as any of these four steps can be done in $O(3^{3k})$ operations, and the modular decomposition contains at most $n$ nodes.
\end{proof}

Although we provide an \FPT\ algorithm, for the modular-width, this algorithm does not directly provide a kernel. Indeed, it has to compute several times the outcome by an exhaustive search on a subgraph. A way to improve this result would be to search for a polynomial kernel for the modular-width.

\subsection{\texorpdfstring{$P_4$}{}-fewness}\label{sec:P4}

Another parameter which can be seen as a distance to cographs is the $P_4$-fewness, which we introduce next. The class of cographs is exactly the class of graphs with no induced $P_4$. Another class of graphs which, in a local sense, contain only a restricted number of $P_4$, introduced by Hoáng \cite{hoangphd}, is the class of {\em $P_4$-sparse} graphs,
which is the class such that no set of five vertices induces more than one $P_4$. Babel and Olariu proposed a generalization of these classes consisting in the class of ``graphs with few $P_4$'s'', or, formally, $(q, q-4)$-graphs~\cite{P4fewness}. A graph is a {\em $(q, t)$-graph} if no set of at most $q$ vertices induces more than $t$ distinct $P_4$’s. With this notation, the class of cographs is the class of $(4,0)$-graphs and the class of $P_4$-sparse graphs is the class of $(5,1)$-graphs.
The {\em $P_4$-fewness} of a graph $G$ is then the smallest integer $q$ such that $G$ is a  $(q,q-4)$-graph~\cite{P4parameterized}. 

The \FPT\ algorithm we will provide for the $P_4$-fewness will rely on the tree decomposition of  $(q, q - 4)$-graphs provided by Jamison and Olariu~\cite{P4fewness,P4fewdecomposition} and reformulated in~\cite{costa2021algorithms} as Primeval Decomposition. We give this decomposition in the following theorem. A {\em spider} is a graph whose vertex set has a partition  $(R, C, S)$, where $C = \{c_1, \ldots , c_m\}$  and $S= \{s_1, \ldots , s_m\}$, for $m \geq 2$, are respectively a clique and a stable set; $s_i$ is adjacent to $c_j$ if and only if $i = j$ (a {\em thin spider}), or $s_i$ is adjacent to $c_j$ if and only if $i \neq j$ (a {\em thick spider}); and every vertex of $R$ is adjacent to each vertex of $C$ and non-adjacent to each vertex of $S$.

\begin{figure}[ht]
    \centering

\centering

\begin{tikzpicture}[scale=0.75]

    \draw (-6,0) ellipse (1cm and 2cm);
    \draw (-2,0) ellipse (1cm and 2cm);
    \draw (2,0) ellipse (1cm and 2cm);


\draw (-6,-1) node[v] (v1){} node[below = .15cm] {$v_1$};
\draw (-6,0) node[] (v7){} node[below = 2.0cm] {$G'$};
\draw (-6,1) node[v] (v2){} node[above = .15cm] {$v_2$};

\draw (-2.0,-1) node[v] (v3){} node[below = .15cm] {$v_3$};
\draw (-2,0) node[] (v8){} node[below = 2.0cm] {$H_1$};
\draw (-2.0,1) node[v] (v4){} node[above = .15cm] {$v_4$};


\draw (2,-1) node[v] (v5){} node[below = .15cm] {$v_5$};
\draw (2,0) node[] (v9){} node[below = 2.0cm] {$H_2$};
\draw (2,1) node[v] (v6){} node[above = .15cm] {$v_6$};

\draw (v1) -- (v3);
\draw (v1) -- (v4);
\draw (v2) -- (v3);
\draw (v2) -- (v4);
%

%

\draw (v3) -- (v5);
\draw (v4) -- (v6);
\draw (v5) -- (v6);

\end{tikzpicture}

\caption{A graph G illustrating Case 5.
}

\label{fig:pseudo spider}
    
\end{figure}
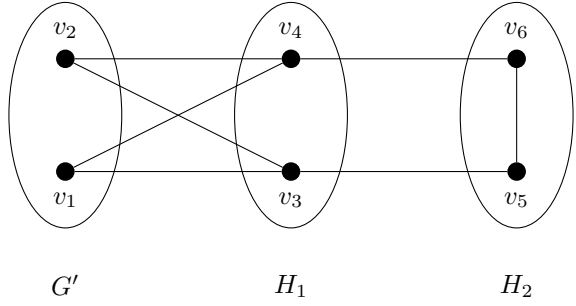

\begin{theorem}[\cite{P4fewness, costa2021algorithms}]\label{P4Decomposition}
Let $G$ be a $(q,q-4)$ graph, then one of the following holds:

\begin{enumerate}
\item $G$ has at most $q$ vertices;
\item $G$ is the disjoint union of two $(q,q-4)$-graphs $G_1$ and $G_2$;
\item $G$ is the join of two $(q,q-4)$-graphs  $G_1$ and $G_2$;
\item $G$ is a spider  $(R,C,S)$  and $G[R]$ is a $(q,q-4)$-graph;
\item $G$ contains a 
non-empty subset of vertices $H$ with $|H|<q$ and bipartition $H=H_1\cup H_2$ such that
$G'=G \setminus H$ is a $(q,q-4)$-graph and every vertex of  $G'$ is adjacent to every vertex of $H_1$ and non-adjacent to every vertex of $H_2$, see Figure~\ref{fig:pseudo spider}. 
\end{enumerate}
\end{theorem}

The decomposition in~\cite{P4fewness, costa2021algorithms} guarantees further properties of $H$ in Case 5 which are not necessary to prove the correctness of our algorithm.

\begin{theorem}\label{thm FPT P4 fewness}
        \MBDG is \FPT\ parameterized by the $P_4$-fewness of the graph.
\end{theorem}

\begin{proof}

Let $G = (V,E)$ be a  $(q,q-4)$-graph.
$G$ can be decomposed by successively applying Theorem~\ref{P4Decomposition}. 
If (2) or (3) holds, we apply the theorem to each components of $G$.
If (4) holds, we apply the theorem to $G[R]$. 
If (5) holds, then we apply the theorem to $G'$.
Finally, (1) corresponds to the leaves of the tree.
This decomposition can be obtained in linear time~\cite{jamison1992tree}.
Then, for each node it is possible to decide the corresponding graph's outcome based on the outcome of its children and the node's corresponding operation. 





\begin{enumerate}
\item If  $|V(G)|\leq q$, we compute the winner by an exhaustive computation of all the possible games. 


\item  If $G = G_1 \cup G_2$,
we compute inductively the outcome of $G_1$ and $G_2$ and
 rely on 
 Lemma~\ref{lemma union}
 to decide the outcome of $G$.


%
\item If $G = G_1 \join G_2$, 
we compute inductively the outcome of $G_1$ and $G_2$ and
 rely on Lemma~\ref{lemma join}  to decide the outcome of $G$.

\item  If 
$G$ is a spider  $(R,C,S)$  and $G[R]$ is a $(q,q-4)$-graph
,
we
 make a further case distinction based on whether the spider is thin or thick.


\begin{enumerate}

\item When the spider is thin,
if Dominator plays first,
Dominator wins,
by playing first in $C$, and then by following a pairing strategy
with pairs $(s_i,c_i)$ for $i \in \{1, \ldots , |C|\}$.
If Staller plays first, 
by iteratively applying Lemma~\ref{unclaimed tail}, where the unclaimed leaves are the vertices of $S$ and their respective private neighbors the vertices of $C$,
we obtain that the outcome of $G$
 is the same as the outcome as $G'$ when Staller plays first. 
Since $G' = G[R]$ is a smaller $(q,q-4)$-graph we 
compute its outcome by induction.

\item When the Spider is thick,
if $|C| \geq 4$,
then $o(G) = \oDom$, with Dominator pairing two disjoint edges in $C$.
Else,
remark that
$G'$ is a module of $G$.
Thus if
$|G'| \geq 3$, we compute inductively the outcome of $G'$, and apply Lemma~\ref{super lemma module edition},
obtaining a new graph $G_2$ with at most nine vertices and same outcome as $G$.
Then we compute $o(G_2)$.
Else, $|G'| < 3$ 
and we compute $o(G)$ as it has at most eight vertices.
\end{enumerate}

\item If
$G$ contains a 
non-empty subset $H$ with bipartition $H_1\cup H_2$, $|H|<q$ and such that
$G'=G \setminus H$ is a $(q,q-4)$-graph  with every vertex of  $G'$ adjacent to every vertex of $H_1$ and non-adjacent to every vertex of $H_2$
, then $G'$ is a module of $G$.
Thus if
$|V(G')| \geq 3$, we compute inductively the outcome of $G'$, and apply Lemma~\ref{super lemma module edition},
obtaining a new graph $G_2$ with at most $q+2$ vertices and same outcome as $G$.
Then we compute $o(G_2)$.
Else, $|V(G')| < 3$ and we compute $o(G)$ as it has at most $q+1$ vertices.
\end{enumerate}

This algorithm runs in time $O\left ( 3^{(q+2)} n \right )$
as any of the steps can be done in $O\left ( 3^{(q+2)} \right )$ operations, and
the primeval decomposition contains at most $n$ nodes, and can be computed
in linear time.
\end{proof}

\section{Distance to cluster}\label{sec:cluster}
We now focus on the distance to cluster, i.e., the minimum number of vertices to delete from our graph to obtain a union of cliques. We first recall that the outcome on a cluster can easily be computed.

\begin{lemma}\label{lemma outcome cluster}
    Let $G$ be a cluster.
\begin{itemize}
    \item If $G$ has no isolated vertices, we have $o(G) = \oDom$.
    \item If $G$ has a single isolated vertex, we have $o(G) = \oN$.
    \item If $G$ has at least two isolated vertices, we have $o(G) = \oS$.
\end{itemize}
\end{lemma}

\begin{proof}
    This result is straightforward since, if $K$ is a clique, it has outcome $\oN$ if it has order $1$, otherwise it has outcome $\oDom$. The result is then given by Lemma~\ref{lemma union}.
\end{proof}

\begin{theorem}\label{theorem FPT dist to cluster}
\MBDG is \FPT, parameterized by the distance to cluster of the graph.
\end{theorem}

\begin{proof}
Let $X$ be a set of vertices such that $G\setminus X$ is a cluster. Let $k=|X|$.
Let $G_1,...,G_s$ be the connected components of $G\setminus X$. All the graphs $G_i$ are cliques.

We first bound the size of the cliques.
First, note that if $G_i$ is a clique, and $v_1, \dots, v_k$ are vertices of $G_i$ with the same neighborhood in $X$, we can, by Lemma~\ref{super lemma module edition} remove all the $v_i$s but two without changing the outcome of the game.
 Therefore, for $1 \le i \le s$, we can suppose $|V(G_i)|\leq 2^{k+1}$, as there are at most $2^k$ different neighborhoods in $X$, and each neighborhood has at most two neighbors in each clique.

The rest of the proof consists in bounding the number of cliques. We define the {\em signature} of a clique $G_i$ as the multiset of the neighborhoods of the vertices of $G_i$ in $X$. Since at most two vertices have the same neighborhood, there are at most $3^{2^{k}}$ possible signatures (for each of the $2^k$ possible neighborhoods, either zero, one or two vertices share this neighborhood). Moreover, if $G_i$ and $G_j$ have the same signature, since we consider multisets, they must have the same number of vertices.
For cliques of size not equal to 2, we prove that we can assume there are at most two other cliques with the same signature. More precisely:

\begin{claim}
Let $G_i$ be a clique of size $\ell$.
\begin{enumerate}
    \item If $\ell=1$, and if at least two other cliques have the same signature, then all of them but two can be removed without changing the outcome of $G$.
    \item If $\ell \geq 3$, and if at least three other cliques have the same signature, then all of them but three can be removed without changing the outcome of $G$.
\end{enumerate}
\end{claim}

\begin{proofclaimitem}
    \begin{enumerate}
        \item Assume $G_i$ has size 1, let $u_i$ be the unique vertex of $G_i$ and $Y$ its neighborhood in $X$. Let $H$ be all the other single vertices having signature $\{Y\}$. Then $H$ is a stable set that is a module of $G$. If $|H|>2$, by Lemma~\ref{super lemma module edition}, it can be replaced by a stable set of size 2. 
        \item Assume now that $G_i$ has size at least 3. Let $\mathcal C$ be all the cliques with the same signature as $G_i$ (with $G_i\in \mathcal C$). Assume that $\mathcal C$ contains at least four cliques. Let $Y$ be all the vertices of $X$ connected to at least one vertex of $G_i$.
        We claim that $\mathcal C$ can be replaced by an edge connected by a join to the set $Y$.
        Let $G'$ be this new graph and $x, y$ be the vertices of the new edge. We prove that $G$ and $G'$ have the same outcome. 
        
        First note that using Lemma~\ref{super lemma}, since $x$ and $y$ have exactly the same neighborhood, we can assume that Dominator will claim $x$ and Staller will claim $y$ in $G'$. In particular, all the set $Y$ is dominated by Dominator.

        Assume first that Staller is winning in $G'$ (going first or second). Then he applies the same strategy in $G$ ignoring the moves of Dominator in $\mathcal C$. Since one can assume in $G'$ that Dominator will claim $x$, the vertex, says $z$, isolated by Staller in $G'$ cannot be a vertex among $Y\cup\{x,y\}$. In particular, $z$ is not connected to any vertex of $\mathcal C$ in $G$ and thus will also be isolated by Staller in $G$. 

         Assume now that Dominator is winning in $G'$ (first or second). As said before, we can already assume that $x$ is claimed by Dominator and $y$ by Staller in $G'$.  Let $\mathcal S'$ be the strategy of Dominator in $G'$ when assuming that $x$ is claimed by Dominator and $y$ by Staller.
        Let $C_1,...,C_m$ be the cliques of $\mathcal C$ ($m\geq 4$). For $j \in \{1,...,m\}$, let $\{u^t_j\}$ with $t\in \{1,....,\ell\}$ be the vertices of $C_j$ such that for a fixed $t$, the vertices $u^t_j$ have the same neighborhood in $X$. 
        While Staller is not playing in $\mathcal C$, Dominator follows her strategy $\mathcal S'$. Without loss of generality, one can assume that $u^1_1$ is the first vertex claimed by Staller in $\mathcal C$. Then Dominator claims $u^1_2$ and thus dominates the whole clique $C_2$. After that, Dominator goes on following $\mathcal S'$ outside $\mathcal C$ and follows a pairing strategy $\mathcal P$ 
        in $\mathcal C$ with the following pairs:
        \begin{enumerate}
            \item For $t\in \{2,...,\ell-2\}$, Dominator adds the pair $(u^t_1,u^t_2)$ in $\mathcal P$.
            \item Dominator adds the pairs $(u^{\ell-1}_2,u^{\ell-1}_3)$ and $(u^{\ell}_2,u^{\ell}_4)$ in $\mathcal P$.
            \item Dominator adds the pair $(u^{\ell-1}_1,u^{\ell}_1)$ in $\mathcal P$.
            \item For $j\in \{3,...,m\}$, Dominator adds the pair $(u^1_j,u^2_j)$ in $\mathcal P$.
        \end{enumerate}

        The pairs from items $(a)$ and $(b)$ ensure that at least one $u^t_i$ is claimed by Dominator for each $t$ while items $(c)$ and $(d)$ ensure that at least one vertex in each clique $C_j$ is claimed by Dominator. Note that the vertices of $\mathcal{C}$ that are not paired can be ignored as they will be dominated by $(c)$ and $(d)$, and cannot dominate new vertices.
        
     This way, Dominator will dominate all the vertices of $\mathcal C$ and will dominate all the vertices of $Y$. All the other vertices of $G$ are present in $G'$ and dominated by another vertex than $x$. Thus, Dominator will dominate the whole graph $G$. \qedclaim
    \end{enumerate}
\end{proofclaimitem}

We now consider the cliques of size 2. We prove that if there are enough cliques of size 2 with the same signature, then we can remove one of them without changing the outcome. By repeating this argument, we prove that we can assume that the number of cliques of size 2 with the same signature is bounded by a function of $k$.

\begin{claim}
    Let $f(k)=(2^{k}+3)3^{2^{k}}+2$. Let $X_1,X_2\subseteq X$. Assume that $G$ has more than $f(k)$ cliques of size 2 with the same signature $(X_1,X_2)$. Let $G'$ be the graph $G$ with one of them removed. Then $G$ and $G'$ have the same outcome.
    \end{claim}

\begin{proofclaimitem}

Let $\mathcal C$ be the set of cliques of size 2 with signature $(X_1,X_2)$. Let $C_1,...C_{|\mathcal C|}$ be the cliques of $\mathcal C$. We assume that $|\mathcal C|>f(k)$. For $j\in \{1,...,\mathcal C\}$, let $u_j$ and $v_j$ be the two vertices of $C_j$ such that all the $u_j$ (respectively all the $v_j$) have neighborhood $X_1$ (resp. $X_2$).  Let $G'$ be the graph where the clique $C_{|\mathcal C|}$ is removed from $G$. We prove that $G$ and $G'$ have the same outcome.
If Dominator wins in $G'$ (first or second), then she wins in $G$ using a pairing strategy on $C_{|\mathcal C|}$ and the same strategy she was using before on the rest of $G$.

Thus assume that Staller wins in $G'$. Let $\mathcal S'$ be a strategy for Staller in $G'$.
We want to prove that Staller wins in $G$. 

We simulate a game in $G$ until $f(k)-1$ cliques of $\mathcal C$ have at least one vertex claimed by a player with Staller that follows $\mathcal S'$ in $G$.
 Without loss of generality, we can assume the $f(k)-1$ cliques played are $C_1$ to $C_{f(k)-1}$.

We observe the following facts:
\begin{enumerate}
    \item If Dominator has claimed one vertex $u_i$ and one vertex $v_j$, then by the Super Lemma, one can assume that for all the empty cliques of $\mathcal C$ there will be one vertex claimed by Staller and one by Dominator, since for all $k \notin \{i,j\}$ and all set $D$ of vertices containing $u_i$ and $v_j$, we have that $D \cup \{u_k\}$ is a dominating set if and only if $D \cup \{v_k\}$ is.
    Then the clique $C_{\mathcal C}$ can be removed, and we obtain the graph $G'$ where Staller will win by hypothesis.
    \item If Dominator is claiming two vertices in $\mathcal C$ with no vertices claimed by Staller in $\mathcal C$ in between, then Dominator has interest to claim two vertices $u_i$ and $v_j$ by Lemma~\ref{lemma: dominated move Maker-Breaker}, and then we are back to Case~1.
    \item Assume that Dominator has already claimed a vertex $v_i$. After this move, if Staller claims a vertex $u_j$ and immediately after Dominator claims $v_j$ then the clique $C_j$ can be removed: Dominator does not dominate any new vertex except $u_j$ and $v_j$ and thus the resulting game is a configuration of $G'$ on which Staller wins. If Dominator claim at some point $v_j$, then it means that $v_j$ had no unclaimed neighbor in $X_2$ nor neighbor claimed by Dominator otherwise it is better for Dominator to claim any neighbor of $v_j$ by Lemma~\ref{lemma: dominated move Maker-Breaker}. In particular, when Staller claims a vertex $u_t$ after that, Dominator has to answer $v_t$ immediately, otherwise Staller isolates a vertex by playing it, and then we can remove the two vertices.
\item Thus one can consider that Dominator has only claimed the vertex $v_1$ in $\mathcal C$ and that Staller claims all the vertices $u_1$ to $u_{f(k)-2}$. When Staller was claiming a vertex $u_j$, Dominator answered outside $\mathcal C$. We will now count the number of optimal moves outside $\mathcal C$. Dominator could claim one of the $k$ vertices of $X$. For the vertices outside $X$, there are at most $3^{2^{k}}$ different signature. For each signature not corresponding to a clique of size 2, there are at most $2^{k}+3$ vertices that are vertices with distinct neighborhood in $X$ or last vertex of a clique that need to be dominated (at most one for each class). For cliques of size 2, once Dominator has claimed the two vertices with distinct neighborhood, from what precedes, she will never claim a vertex in a clique except if Staller is claiming a vertex and threats to isolate a vertex. Thus, when Staller is claiming a vertex of $\mathcal C$, Dominator will play at most twice in each class of cliques of size 2. Therefore, in total there are at most $(2^{k}+3)3^{2^{k}}=f(k)-2$ interesting vertices. Once all these vertices have been claimed by Dominator, she should answer in $\mathcal C$ directly after the next move of Staller. Since there is still a vertex $u_{f(k)-1}$ unclaimed she can claim it, and we are back to Case~1. \qedclaim
\end{enumerate}
\end{proofclaimitem}

To conclude the proof, we proved that, for each signature, we can reduce the graph $G$ to a graph in which at most $f(k)$ cliques share a signature. As there are at most $3^{2^k}$ signatures in total, and as each clique has at most $2^{k+1}$ vertices, we can reduce the graph to a kernel having $f(k)3^{2^k}2^{k+1} = O(3^{2^{k+2}})$ vertices.

We recall that having a kernelization algorithm is equivalent to be \FPT. Therefore, \MBDG is \FPT\ parameterized by the distance to cluster.
\end{proof}

\section{Feedback edge number}\label{sec:fes}
We now focus on the  feedback edge number. A set of edges $S$ is a {\em feedback edge set} if $G \setminus S$ is a forest. The {\em feedback edge number} of a graph $G$, denoted by $fen(G)$, is then defined as the size of a minimum feedback edge set. 

When Staller starts, we will prove that \MBDG has a linear kernel.
Up to considering the $|V|$ possible first moves, we still obtain an \FPT\ algorithm in the general case. 

\begin{theorem}\label{thm: MB Dom fes}
\MBDG has a linear kernel, parameterized by the feedback edge number, after the first move of Dominator. 
\end{theorem}

\begin{corollary}\label{cor: MB Dom fes}
    \MBDG is \FPT\ parameterized by the feedback edge number of the graph.
\end{corollary}

The proof mostly relies on the fact that, if a graph has a bounded feedback edge number, either it has a lot of leaves, or most of its vertices belong to long induced paths. More formally, we prove the following structural lemma:

\begin{lemma}\label{vertices high degree fes}
    Let $G$ be a graph, with at most one leaf. Let $L$ be the set of vertices of $G$ of degree $3$ or more. We have $|L| \le 2 fen(G)-1$.
\end{lemma}

\begin{proof}
    Without loss of generality, we can suppose $G$ connected. Let $H$ be the graph obtained by contracting any vertex of degree $2$ in $G$ and by removing its leaf (if any). Note that $fen(H) = fen(G)$ and that the vertices of $H$ are exactly $L$ or $L\setminus \{x\}$ if $x$ is the neighbor of the removed leaf. All the vertices of $H$ have degree at least $3$, so we have $|E(H)| \ge \frac{3}{2} |V(H)|$. Since $|E(H)| - |V(H)| +1 = fen(H)$, we have by subtracting this equality to the previous one $|V(H)| - 1 \ge \frac{3}{2} |V(H)| - fen(H)$. Finally, we obtain $2fen(H) - 2 = 2fen(G) - 2 \ge |V(H)|$. As putting back the leaf can at most create one vertex of degree $3$, we have $2fen(G) -1 \ge |L|$.
\end{proof}

The main lemma of the proof is the following one, which states that long paths of the graph can be shortened. Its proof relies on the fact that after the first move in a path, the path can be split, then pending paths can be shortened and moves will be made such that the isolated paths are winning for Dominator.  

\begin{figure}[ht]
    \centering

\begin{tikzpicture}

    \draw (0,0) circle (1) ;
    \draw (0,0) node[inv] (G){\large $G$} ;

    \foreach \i/\label in {1/1, 2/2, 3/3, 4/4, 5/5, 6/6, 7/7, 8/8, 9/9} {
    \coordinate (V\i) at ({360/16 * (\i - 1)}:2);
\ifthenelse{\i<6}{\draw (V\i) node[v] (1.5pt){} node[above right] {$v_{\label}$};}{\draw (V\i) node[v] (1.5pt){} node[above left] {$v_{\label}$};}
    }

\draw (-1,0) node[v] (u){} node[right=.15cm] {$u$};
\draw (1,0) node[v] (v){} node[left=.15cm] {$v$};

    \foreach \i/\j in {1/2, 2/3, 3/4, 4/5, 5/6, 6/7, 7/8, 8/9} {
    \draw (V\i) -- (V\j);
    }
\draw (V9) -- (u);
\draw (V1) -- (v);

    \draw (7,0) circle (1) ;
    \draw (7,0) node[inv] (G){\large $G$} ;

    \foreach \i/\label in {1/1, 2/2, 3/3, 4/4, 5/5, 6/6, 7/7} {
    \coordinate (V\i) at ({360/12 * (\i - 1)}:2) ;
\ifthenelse{\i<5}{\draw (V\i)+(7,0) node[v](u\i) {} node[above right] {$v_{\label}$};}{\draw (V\i)+(7,0) node[v](u\i) {} node[above left] {$v_{\label}$};}
    }

\draw (6,0) node[v] (u){} node[right=.15cm] {$u$};
\draw (8,0) node[v] (v){} node[left=.15cm] {$v$};

    \foreach \i/\j in {1/2, 2/3, 3/4, 4/5, 5/6, 6/7} {
    \draw (u\i) -- (u\j);
    }

\draw (V7)+(7,0) -- (u);
\draw (V1)+(7,0) -- (v);

\end{tikzpicture}
    
    \caption{Two equivalent graphs, using Lemma~\ref{shorten any path in graph}. On the left, it is $G^{u,v}_9$. On the right, it is $G^{u,v}_7$}
    \label{fig:shorten path}
\end{figure}
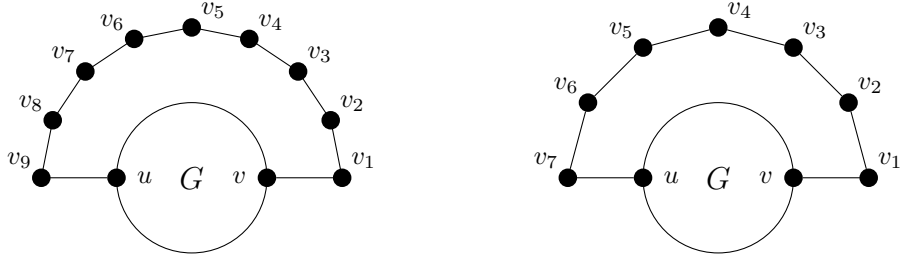

\begin{lemma}\label{shorten any path in graph}
    Let $G$ be a graph, and $u,v$ be two vertices of $G$. We denote by $G^{u,v}_k$ the graph obtained from $G$ by adding a path on $k$ vertices to $G$ and connecting $u$ and $v$ to the two extremities of the path. If $k \ge 7$, then $G^{u,v}_k$ and $G^{u,v}_{k+2}$ have the same outcome, even if some vertices of $G$ but not on the path, have already been played by Dominator. See Figure~\ref{fig:shorten path} for an illustration.
\end{lemma}

Since isolated paths will be created by making use of Lemma~\ref{removing vertices Dominator},
we first show that paths on which an extremity is claimed by Dominator are winning for Dominator.

\begin{lemma}\label{lem: path extremity dom}
    Let $P_n$ be a path of order $n$. We have $o(P_n, \{v_1\}, \emptyset) = o(P_n, \{v_1, v_n\}, \emptyset) = \oDom$ where $v_1$ and $v_n$ are the extremities of $P_n$.
\end{lemma}

\begin{proof}
It is sufficient to prove this lemma with only $v_1$ claimed by Dominator, since when Dominator has also claimed $v_n$, the position is more favorable for Dominator then the one where she has not according to Observation~\ref{obs: free vertex}. 

    Let $P_n$ be a path of order $n$. Denote its vertices $v_1, \dots, v_n$ such that $(v_i,v_{i+1}) \in E(P_n)$. Suppose that Dominator has already claimed $v_1$. If $n$ is odd, $\mathcal P = \{ (v_{2i}, v_{2i+1})\}_{1 \le i \le \frac{n-1}{2}}$ 
    provides a winning pairing strategy for Dominator.
    If $n$ is even, $\mathcal{P} = \{ (v_{2i-1}, v_{2i})\}_{2 \le i \le \frac{n}{2}}$ 
        provides a winning pairing strategy for Dominator.
Thus, $P_n$ has outcome $\oDom$.
\end{proof}

The proof of Lemma~\ref{shorten any path in graph} uses another more specific lemma in order to deal with long pending paths:

\begin{lemma}\label{shorten path Dominator}
Let $G$ be a graph and $k\geq 3$. We denote by $G^u_k$ the graph obtained by adding a pending path on $k$ vertices $(v_1, \dots, v_k)$ to $G$ and connecting $v_1$ to $u$. Let $G_k=(G^u_k,D,S)$ be a position on $G^u_k$ such that the only vertex claimed on the vertices $\{v_1,...,v_k\}$ is $v_k$ and is claimed by Dominator.
If Dominator wins in the position $G_k$, then she wins in all the positions $G_{k'}=(G^u_{k'},D\setminus \{v_k\}\cup \{v_{k'}\},S)$, for any $k'\geq 1$.
\end{lemma}

\begin{proof}
When $k'<k$ it is clear : the minimal dominating sets for Dominator in $G_k$ are the minimal dominating sets of $G_{k'}$ (when restricted to the unclaimed vertices). Thus, Dominator can just follow her strategy in the smaller graph. If at some point she has to claim a vertex not in $G_{k'}$, she just claims any unclaimed vertex.

We now prove the result for $k'=k+1$. By induction, it will prove the result for all $k'\geq k$.

Dominator follows her strategy in $G_k$ until an element of $X=\{v_1,...,v_k\}$ is claimed. Note that up to change the sets $D$ and $S$, we can assume that it is the first move.
\begin{itemize}
    \item Assume first that it is Dominator that claims a vertex $v_j$ in $X$. Then, using Lemma~\ref{removing vertices Dominator} the game $G_k$ (resp. $G_{k+1}$) is split into two games: $G_j$ and a path with $k-j+1$ (resp. $k-j+2$) vertices with only the two end vertices claimed by Dominator. Since playing $v_j$ was an optimal move for Dominator in $G_k$, she is winning playing second in the game $G_k$ with $v_j$ claimed by Dominator. Therefore, she is winning second in $G_j$ (by Lemma~\ref{lemma union}). Since she is also winning second in any path with extremities claimed by Dominator by Lemma~\ref{lem: path extremity dom}, she is winning in $G_{k+1}$ by claiming $v_j$.

    \item Assume now that it is Staller that claims a vertex $v_j$ in $X$. If $j>1$, then Dominator claims $v_{j-1}$ (that is unclaimed by hypothesis). Then as before the game $G_{k+1}$ is split into a path with one extremity claimed by Dominator (where Dominator wins second, by Lemma~\ref{lem: path extremity dom}) and the game $G_{j-1}$ that is winning for Dominator playing second (since $1\leq j-1\leq k$ and the first remark of the proof). Thus, the union is winning by Dominator playing second by Lemma~\ref{lemma union}.

   Assume now $j = 1$: Staller claims $v_1$. If $u$ is unclaimed, Dominator claims $u$ and wins since $(G,D\cup \{u\},S)$ is winning for Dominator playing second, as a subgraph of $G_k$ which had outcome $\oDom$ and a path with one extremity claimed by Dominator, which has outcome $\oD$ by Lemma~\ref{lem: path extremity dom}. If $u$ was already claimed by Dominator, Dominator can answer in $G$, and we have the same result.
    Otherwise, if it is claimed by Staller, then in both games, Dominator has to answer in $v_2$ since $k\geq3$ and $v_1$ need to be dominated. Then both games are  split into two components 
    whose outcome is  $\oDom$
    by hypothesis for $G$ and Lemma~\ref{lem: path extremity dom}.
\end{itemize}

\end{proof}

We can now prove Lemma~\ref{shorten any path in graph}.

\begin{proof}[Proof of Lemma~\ref{shorten any path in graph}]
Denote by $P_k$ the path attached to $u$ and $v$, and $D$ the set of vertices already claimed by Dominator in $G$.

    Assume first that Dominator has a winning strategy $\strat$ in $G^{u,v}_k$ where she has already claimed $D$.  We prove that Dominator wins in $G^{u,v}_{k+2}$. For her first move, if it is her turn, and while Staller plays in $G$, Dominator answers following $\strat$ in $G$.  If before any vertex is played on $P_{k+2}$, Dominator played $u$ or $v$, by Lemma~\ref{removing vertices Dominator}, we can split this vertex. Then, by applying Lemma~\ref{shorten path Dominator}, we can remove two vertices on this path, and we obtain the result directly. Thus, we can consider that $u$ and $v$ are either unclaimed or claimed by Staller. We consider the first time a vertex has to be played on $P_k$ or is played on $P_{k+2}$.
\begin{itemize}
    \item If $\strat$ wants Dominator to claim in $P_k$, we denote the vertices of $P_k$ by $v_1, \dots, v_k$, in such a way that $\strat$ makes Dominator plays $v_i$ with $i$ minimal. Dominator plays $v_i$ in $P_{k+2}$. Now, by Lemma~\ref{removing vertices Dominator}, we can split the path into two, one of length $i$, and one of length $k+3-i$, by removing $v_i$ and replacing it with two leaves claimed by Dominator and connected to its neighbors. Then, by Lemma~\ref{shorten path Dominator}, as either $i$ or $k+3-i$ is greater than $3$, (since $k \ge 7$), we can shorten one of these two paths by two vertices (or both by one). Then, using again Lemma~\ref{removing vertices Dominator}, the resulting graph has the same outcome as if $v_i$ was played in $P_k$, which concludes the proof.
    \item If Staller claims first in $P_{k+2}$, we denote the vertices of $P_{k+2}$ by $v_1, \dots, v_{k+2}$, in such a way that Staller played $v_i$ with $i$ minimal. Suppose $v_1$ is connected to $u$. The following happens depending on the value of $i$: \begin{itemize}
        \item If $i\ge 2$, Dominator, in $\strat$ has to play a neighbor of $v_i$, otherwise Staller claims either $v_{i-1}$ or $v_{i+1}$ and threaten to isolate a vertex with his next move. 
        
        Thus, Dominator can claim the same vertex in $P_k$ (either $v_{i-1}$ or $v_{i+1}$). Now, by Lemma~\ref{removing vertices staller}, $v_i$ can be removed. Since $k \ge 7$, On the path $v_{i+1}, \dots, v_{k+2}$, we can apply either Lemma~\ref{shorten path Dominator}, if Dominator has claimed $v_{i+1}$, or Lemma~\ref{unclaimed tail}, if she has played $v_{i-1}$, which gives the result by induction.
        \item If $i = 1$, Dominator has to claim a vertex in $\{u, v_2, v_3\}$, otherwise Staller can claim $v_2$ and threaten to isolate $v_1$ by claiming $u$ or $v_2$ by claiming $v_3$ (or wins if he already has claimed $u$). Thus, after this move, we can split this vertex using Lemma~\ref{removing vertices Dominator}. Once again, if Dominator has claimed $u$, we conclude by Lemma~\ref{unclaimed tail}. If Dominator has claimed $v_2$ or $v_3$, we conclude by Lemma~\ref{shorten path Dominator}. 
    \end{itemize}
\end{itemize}

Suppose now that Dominator has a winning strategy in $G_{k+2}^{u,v}$ where she has already claimed the vertices of $D$. We prove that Dominator wins in $G^{u,v}_{k}$. For her first move, if it is her turn, and while Staller plays in $G$, Dominator answers following $\strat$ in $G$.  If before any vertex is played on $P_{k}$, Dominator played $u$ or $v$, by Lemma~\ref{removing vertices Dominator}, we can split this vertex. Then, by applying Lemma~\ref{shorten path Dominator}, applied on $G_{k}^{u,v}$, we can add two vertices on this path, and we obtain the result directly. Thus, we can consider that $u$ and $v$ are either unclaimed or claimed by Staller. We consider the first time a vertex has to be played on $P_{k+2}$ or is played on $P_{k}$.
\begin{itemize}
    \item If $\strat$ wants Dominator to claim in $P_{k+2}$, we denote the vertices of $P_{k+2}$ by $v_1, \dots, v_{k+2}$, in such a way that $\strat$ makes Dominator plays $v_i$ with $i$ minimal. Dominator plays $v_i$ in $P_{k}$. Now, by Lemma~\ref{removing vertices Dominator}, we can split the path into two, one of length $i$, and one of length $k+1-i$, by removing $v_i$ and replacing it with two leaves claimed by Dominator and connected to its neighbors. Then, by Lemma~\ref{shorten path Dominator}, as either $i$ or $k+1-i$ is greater than $3$, (since $k \ge 7$), we can add two vertices to one of these two paths. Then, using again Lemma~\ref{removing vertices Dominator}, the resulting graph has the same outcome as if $v_i$ was played in $P_{k+2}$, which concludes the proof.
    \item If Staller claims first in $P_{k}$, we denote the vertices of $P_{k}$ by $v_1, \dots, v_{k}$, in such a way that Staller played $v_i$ with $i$ minimal. Suppose $v_1$ is connected to $u$. The following happens depending on the value of $i$: \begin{itemize}
        \item If $i\ge 2$, Dominator, in $\strat$ has to play a neighbor of $v_i$, otherwise Staller claims either $v_{i-1}$ or $v_{i+1}$ and threaten to isolate a vertex with his next move. 
        
        Thus, Dominator can claim the same vertex in $P_k$ (either $v_{i-1}$ or $v_{i+1}$). Now, by Lemma~\ref{removing vertices staller}, $v_i$ can be removed. Since $k \ge 7$, On the path $v_{i+1}, \dots, v_{k}$, we can apply either Lemma~\ref{shorten path Dominator}, if Dominator has claimed $v_{i+1}$, or Lemma~\ref{unclaimed tail}, if she has played $v_{i-1}$, which gives the result by induction.
        \item If $i = 1$, Dominator has to claim a vertex in $\{u, v_2, v_3\}$, otherwise Staller can claim $v_2$ and threaten to isolate $v_1$ by claiming $u$ or $v_2$ by claiming $v_3$ (or wins if he already has claimed $u$). Thus, after this move, we can split this vertex using Lemma~\ref{removing vertices Dominator}. Once again, if Dominator has claimed $u$, we conclude by Lemma~\ref{unclaimed tail}. If Dominator has claimed $v_2$ or $v_3$, we conclude by Lemma~\ref{shorten path Dominator}. 
    \end{itemize}
\end{itemize}

Finally, we proved that if Dominator has a winning strategy on one of $G_k^{u,v}$ or $G_{k+2}^{u,v}$, she has one on both, which proves that they have the same outcome.
\end{proof}

Using these lemmas, we can finally prove Theorem~\ref{thm: MB Dom fes}.

\begin{proof}[Proof of Theorem~\ref{thm: MB Dom fes}  and Corollary~\ref{cor: MB Dom fes}]
    Let $G$ be a graph and let $k = fen(G)$. Let $n = |V(G)|$. Up to considering the $O(n)$ positions reached after the first move of Dominator, suppose that it is Staller's turn. By Lemma~\ref{unclaimed tail}, Staller claims all the unclaimed neighbors of leaves forcing answers to Dominator after each move on the attached leaf. By applying Lemma~\ref{removing vertices staller}, and by removing the leaves connected to the first move of Dominator, which are now in no minimal dominating sets, we can assume that the resulting graph has at most one leaf, and if so, it has been played by Dominator. 

    The resulting graph, by Lemma~\ref{vertices high degree fes}, has at most $2k-1$ vertices of degree $3$ or more. Now, by applying Lemma~\ref{removing vertices Dominator}, we can split the first vertex of Dominator to transform it into a set of leaves, without changing the degree of any other vertex, nor the number of unclaimed vertices. Now, the only vertices that are already claimed by Dominator are leaves of pending paths. Thus, by Lemma~\ref{shorten any path in graph} and Lemma~\ref{shorten path Dominator}, we can shorten any path containing only vertices of degree $2$ in this graph to a path of order at most $7$. 
    
   The resulting graph still has a feedback edge set of order $k$ and at most one leaf. Therefore, after removing the feedback edge set, we obtain a forest with at most $2k$ vertices of degree at least $3$, and at most $2k$ leaves. Since paths of vertices of degree $2$ have at most length $7$, this tree has at most $7(4k-1)$ vertices of degree two. Finally, it has $O(k)$ vertices. So, the given kernel is linear.
    
    Finally, we proved that the game has a linear kernel after the first move of Dominator. Up to consider the different first moves for Dominator, which only adds a factor $n$ to the complexity, there is an \FPT{} algorithm to determine the winner in graphs of having a feedback edge set of size $k$.
\end{proof}

\section{Conclusion and open questions}

In this paper we have studied the parameterized complexity of the Maker-Breaker Domination game from the point of view of various parameters. For the natural parameter of the number of moves we have proved that determining whether Staller can isolate a vertex in $k$ moves and determining whether Dominator can claim a dominating set in $k$ moves are respectively  \W[1]-complete and \W[2]-complete parameterized by $k$. We have actually proved a stronger result: in a Maker-Breaker positional game, determining whether Breaker can claim a transversal in $k$ moves is \W[2]-complete. This contrasts with the case when Maker has to claim a hyperedge in $k$ moves, which was known to be \W[1]-complete.

Moreover, since {\shortDD} and {\shortDS} can be expressed as first order logic formula, they can be solved in \FPT~time, parameterized by the number of moves $k$, on nowhere dense graphs due to the meta-theorem from Grohe \etal~\cite{Grohe2013}. 
Similarly, using Courcelle's meta-theorem~\cite{Courcelle1990}, we obtain that \shortDD and \shortDS can be computed in time $O(f(k,t)\cdot n)$ on graphs of treewidth~$t$.
However, this approach requires the number of moves to be taken as a parameter, and the complexity of \MBDG only parameterized by the treewidth of the input graph is still open. Even on treewidth~2 graphs, as mentioned in~\cite{Intervalle}, it is still open to know if we have a polynomial algorithm to compute the outcome (but a polynomial-time algorithm is known for outerplanar graphs). Note that it is unlikely to obtain a positive result using model checking arguments since QBF is already \PSPACE-complete on bounded treewidth formulas~\cite{QBFwidth}.

We have also studied the parameterized complexity with regard to several structural parameters, among which the modular width, the $P_4$-fewness, the distance to cluster, and the feedback edge number, for which we provide \FPT\ algorithms.
Concerning the modular width, the proof provided here does not provide us a kernel. Thus, one could try to obtain a polynomial kernel for this parameter. Since this method mostly relies on an application of Lemma~\ref{super lemma}, and since to the best of our knowledge, replacing modules in a graph with a smaller graph with the same outcome was never used in the complexity study of positional games, 
this could be a new approach for proving that a positional game is \FPT\ parameterized by the modular-width, by providing lemmas similar to Lemma~\ref{super lemma module edition} for other positional games.

We studied the distance to cluster, as clusters are graphs in which Dominator wins easily, we can wonder what happens to other classes of graphs where the winner can be computed in polynomial time. For instance, considering the feedback vertex set instead of the feedback edge set would give a distance to forests in terms of number of vertices and not edges. This parameter seems a bit difficult to handle, since we cannot use it to bound the number of vertices of large degree, as we did for the feedback edge set. A first step could be the distance to union of stars.

Finally, in~\cite{Oijid2025}, it is proved that the Domination game remains \PSPACE-complete in bounded degree graphs, and in~\cite{Intervalle}, that Dominator always wins in regular graphs. Considering the complexity of \MBDG played on subcubic graphs would be an interesting question since moves of Dominator can quickly be forced. Moreover, the underlying hypergraph of the closed neighborhoods would have rank~$4$, for which it is \PSPACE-complete in general to compute the outcome of a positional game~\cite{galliot2025}.


\bibliographystyle{plain}
\bibliography{journal}

\end{document}